\documentclass[10pt]{article}
\usepackage{tikz}
\usetikzlibrary{calc}
\usepackage{amsmath, upgreek}
\usepackage{amssymb}
\usepackage{color}
\usepackage{amscd}
\usepackage{xspace}
\usepackage{verbatim}
\usepackage{graphicx}
\usepackage{color}
\renewcommand{\theequation}{\thesection.\arabic{equation}
}
\newcommand{\mysection}[1]{
\section{#1}\setcounter{equation}{0}}
 \usepackage{cite}

\begin{document}


\newcommand{\txt}[1]{\;\text{ #1 }\;}
\newcommand{\tbf}{\textbf}
\newcommand{\tit}{\textit}
\newcommand{\tsc}{\textsc}
\newcommand{\trm}{\textrm}
\newcommand{\mbf}{\mathbf}
\newcommand{\mrm}{\mathrm}
\newcommand{\bsym}{\boldsymbol}
\newcommand{\scs}{\scriptstyle}
\newcommand{\sss}{\scriptscriptstyle}
\newcommand{\txts}{\textstyle}
\newcommand{\dsps}{\displaystyle}
\newcommand{\fnz}{\footnotesize}
\newcommand{\scz}{\scriptsize}
\newcommand{\be}{\begin{equation}}
\newcommand{\bel}[1]{\begin{equation}\label{#1}}
\newcommand{\ee}{\end{equation}}
\newcommand{\eqnl}[2]{\begin{equation}\label{#1}{#2}\end{equation}}
\newcommand{\barr}{\begin{eqnarray}}
\newcommand{\earr}{\end{eqnarray}}
\newcommand{\bars}{\begin{eqnarray*}}
\newcommand{\ears}{\end{eqnarray*}}
\newcommand{\nnu}{\nonumber \\}
\newtheorem{subn}{\name}
\renewcommand{\thesubn}{}
\newcommand{\bsn}[1]{\def\name{#1}\begin{subn}}
\newcommand{\esn}{\end{subn}}
\newtheorem{sub}{\name}[section]
\newcommand{\dn}[1]{\def\name{#1}}   
\newcommand{\bs}{\begin{sub}}
\newcommand{\es}{\end{sub}}
\newcommand{\bsl}[1]{\begin{sub}\label{#1}}
\newcommand{\bth}[1]{\def\name{Theorem}
\begin{sub}\label{t:#1}}
\newcommand{\blemma}[1]{\def\name{Lemma}
\begin{sub}\label{l:#1}}
\newcommand{\bcor}[1]{\def\name{Corollary}
\begin{sub}\label{c:#1}}
\newcommand{\bdef}[1]{\def\name{Definition}
\begin{sub}\label{d:#1}}
\newcommand{\bprop}[1]{\def\name{Proposition}
\begin{sub}\label{p:#1}}

\newcommand{\aand}{\quad\mbox{and}\quad}
\newcommand{\M}{{\cal M}}
\newcommand{\A}{{\cal A}}
\newcommand{\B}{{\cal B}}
\newcommand{\I}{{\cal I}}
\newcommand{\J}{{\cal J}}
\newcommand{\D}{\displaystyle}
\newcommand{\RR}{ I\!\!R}
\newcommand{\C}{\mathbb{C}}
\newcommand{\cC}{\mathcal{C}}
\newcommand{\cP}{\mathcal{P}}
\newcommand{\R}{\mathbb{R}}
\newcommand{\Z}{\mathbb{Z}}
\newcommand{\N}{\mathbb{N}}
\newcommand{\T}{{\rm T}^n}
\newcommand{\cuad}{{\sqcap\kern-.68em\sqcup}}
\newcommand{\abs}[1]{\mid #1 \mid}
\newcommand{\norm}[1]{\|#1\|}
\newcommand{\equ}[1]{(\ref{#1})}
\newcommand\rn{\mathbb{R}^N}
\renewcommand{\theequation}{\thesection.\arabic{equation}}
\newtheorem{definition}{Definition}[section]
\newtheorem{theorem}{Theorem}[section]
\newtheorem{proposition}{Proposition}[section]
\newtheorem{example}{Example}[section]
\newtheorem{proof}{proof}[section]
\newtheorem{lemma}{Lemma}[section]
\newtheorem{corollary}{Corollary}[section]
\newtheorem{remark}{Remark}[section]
\newcommand{\bremark}{\begin{remark} \em}
\newcommand{\eremark}{\end{remark} }
\newtheorem{claim}{Claim}

\newcommand{\cH}{\mathcal{H}}
\newcommand{\cL}{\mathcal{L}}
\newcommand{\cO}{\mathcal{O}}
\newcommand{\cA}{\mathcal{A}}
\newcommand{\cJ}{\mathcal{J}}
\newcommand{\cF}{\mathcal{F}}
\newcommand{\cE}{\mathcal{E}}
\newcommand{\BBE}{\mathbb{E}}

 \begin{center}
 \large   Existence of solutions to  elliptic equations  involving   regional \\[2mm]
 fractional Laplacian with order   $(0,\frac12]$
 \end{center}

 \medskip

\centerline{\small    Huyuan Chen\footnote{chenhuyuan@yeah.net},\ \ \quad
Huihuan Peng\footnote{penghuihuan00@126.com},\ \ \quad
Yanqing Sun\footnote{sunyanqing0114@126.com} }
\bigskip
{

  \centerline{\small $^1$ Shanghai Institute for Mathematics and Interdisciplinary Sciences,\
Fudan University, }
  \centerline{\small  Shanghai ‌200433, PR  China}

\medskip

  \centerline{\small $^{2,3}$ Department of Mathematics, Jiangxi Normal University, Nanchang,}
  \centerline{\small   Jiangxi 330022, PR China}

}
\medskip

\begin{abstract}
 Our purpose of this paper is to investigate  positive solutions of the
elliptic equation with regional fractional Laplacian
$$
    ( - \Delta )_{B_1}^s u  +u= h(x,u)  \quad
     {\rm  in} \ \,   B_1,\qquad   u\in C_0(B_1),
    $$
  where     $( - \Delta )_{B_1}^s$ with  $s\in(0,\frac12]$ is the regional fractional Laplacian and $h$ is the nonlinearity.

 Ordinarily, positive solutions vanishing at the boundary are not anticipated to be derived for the equations with regional fractional Laplacian of order $s\in(0,\frac12]$. Positive solutions are obtained when the nonlinearity assumes the following two models:
  $h(x,t)=f(x)$
or  $h(x,t)=h_1(x)\, t^p+ \epsilon h_2(x)$,  where $p>1$, $\epsilon>0$ small  and  $f, h_1, h_2$ are H\"older continuous, radially symmetric and decreasing functions under suitable conditions.

\end{abstract}
 \medskip


\vspace{1mm}
\hspace{.05in}
\smallskip

\setcounter{equation}{0}
\section{Introduction}

Let  $s\in (0,1)$, $\Omega$ be an  $C^2$ domain in $\R^N$ with $N\geq2$,
$ (-\Delta)^s_\Omega$ be the regional fractional Laplacian  defined by
$$ (-\Delta)^s_\Omega  u(x) =c_{N,s} \lim_{\varepsilon\to0^+} \int_{\Omega\setminus B_\varepsilon(x)}\frac{ u(x)-u(z)}{|z-x|^{N+2s}}\,dz, $$
where $B_r(x)$ is the ball with radius $r$ and the center at $x$, particularly, denote $B_r=B_r(0)$,
here $c_{N,s}>0$ is the normalized constant of  fractional Laplacian $(-\Delta)^s_{\R^N}$ (simply we use the notation $(-\Delta)^s$), see \cite{EGE}.  

 In recent years, nonlocal problems have been increasingly studied across various fields such as physics models, operations research, queuing theory, mathematical finance, and risk estimation (see \cite{CSS}). The regional fractional Laplacian is a representative operator associated with the generator of a censored stable process. From a probabilistic perspective, a symmetric $2s$-stable process in $\mathbb{R}^N$ that is killed upon exiting a domain $\Omega$ is referred to as a symmetric $2s$-stable process confined to $\Omega$. Bogdan, Burdzy, and Chen \cite{BBC} (see also Guan and Ma \cite{CKS, GM}) extended this class of processes to construct a version of a strong Markov process, known as the censored symmetric stable process. A censored stable process in an open set $\Omega \subset \mathbb{R}^N$ is obtained by suppressing the jumps of a symmetric stable process from $\Omega$ to $\mathbb{R}^N \setminus \Omega$. It is worth noting that censored stable processes exhibit distinctive properties that highlight differences between the cases $s \in (\frac{1}{2}, 1)$ and $s \in (0, \frac{1}{2}]$ (see \cite[Theorem 1]{BBC}):

(i) for $s\in(\frac12,\, 1)$,  the censored symmetric $2s$-stable process in $\Omega$ has a finite lifetime and will approach $\partial \Omega$;

(ii) for  $s\in (0, \, \frac12]$, the censored symmetric $2s$-stable process in $\Omega$   is conservative and will never approach  $\partial \Omega$.

On the analysis side, interesting new phenomena occur in relation to
elliptic problems involving the regional fractional Laplacian.
   Let
  ${\bf H}^s_0(\Omega)$ be the closure of  $C^\infty_c(\Omega)$  under the semi-norm that
  $$\|u\|_{s,\Omega} =\sqrt{\int_{\Omega\times \Omega} \frac{(u(x)-u(y))^2}{|x-y|^{N+2s}}dxdy}.$$
The authors in \cite{M} showed that {\it
   for $s\in(\frac12,\, 1)$,  Hilbert space ${\bf H}^s_0(\Omega)$ has    zero trace, 
 while for  $s\in (0, \, \frac12]$, Hilbert space ${\bf H}^s_0(\Omega)$  has no   zero trace. }

 Let $H^s_0(\Omega)$ be the closure of  $C^\infty_c(\Omega)$  under the norm that
  $$\|u\|_{s,\Omega} =\sqrt{\int_{\Omega\times \Omega} \frac{(u(x)-u(y))^2}{|x-y|^{N+2s}}dxdy+\int_{\Omega}u^2dx }. $$
It is worth noting that  function $1\in H^s_0(\Omega)$,  $H^s_0(\Omega)={\bf H}^s_0(\Omega)\cap L^2(\Omega)$ and it also has no    zero trace for $s\in(0,\frac12]$. 
  This means  it is delicate to determine, for $s\in(0,\frac12]$,  whether there is  a nontrivial solution of the related Schr\"odinger equation with regional fractional Laplacian even in a ball
\begin{equation}\label{eq 1.1}
    \left\{ \arraycolsep=1pt
    \begin{array}{lll}
    ( - \Delta )_{B_1}^s u+u=h(x,u) \quad \
    &{\rm in} \ \, B_1,\\[2mm]
    \qquad \ \,    u\in C_0(B_1),
    \end{array}
    \right.
    \end{equation}
where     $h: B_1\times \R\to \R$ is a measurable function and {\it $C_0( B_1)$ is the set of  function continuous in $\bar B_1$,
 which vanishes at the boundary $\partial B_1$.} Our primary objective in this paper is to investigate the existence of nontrivial solutions to (\ref{eq 1.1}) when the nonlinearity $ h $ takes typical models, such as non-homogeneous terms.

 When $s\in(\frac12,\, 1)$,  \cite{CK,CKS} provide estimates on the heat kernel and Green kernel related to the regional fractional Laplacian,  \cite{G} builds a formula of  integration by part for regional fractional  Laplacian,  \cite{C} extends this formula to solve regional fractional  problem with  inhomogeneous terms.
Via building the formula of integral by part  and related  embedding results,   \cite{C} obtains the existence of  solutions to
\begin{equation}\label{eq 1.1-y}
\arraycolsep=1pt\left\{
\begin{array}{lll}
 \displaystyle  (-\Delta)^s_\Omega   u=f\quad & {\rm in}\ \,  \Omega,\\[1.5mm]
\phantom{ (-\Delta)^\alpha  }
 \displaystyle   u=0\quad & {\rm on}\   \partial\Omega
 \end{array}\right.
\end{equation}
 for $s\in(\frac12,1)$ when $\Omega$ is a bounded regular domain.  For further study of regional fractional Laplacian with $s\in (\frac12,1)$,    refer to  \cite{CH,AGV}  for boundary blowing-up solutions,   \cite{F} for boundary regularity, \cite{AFT} for related Hopf Lemma, \cite{FT} for existence of weak solutions with critical semilinear term.   More related study see      \cite{AT,DFW,BGU} and references therein on related topics  involving the regional fractional   Laplacian.

 \smallskip

 For $s\in(0,\frac12]$, the structure of solutions of elliptic equations are very challenging.
 The authors in\cite{CW} showed the nonexistence of solutions to Poisson problem
$$(-\Delta)^s_\Omega u=1\quad{\rm in}\ \, \Omega$$
 and nonexistence of positive solutions to Lane-Emden  equation
  \begin{equation}\label{eq 1.2}
  (-\Delta)^s_\Omega u=u^p\quad{\rm in}\ \, \Omega.
  \end{equation}
Under the assumption of  $\int_\Omega f dx=0$,
   \cite{T} shows the existence of weak solutions to Poisson problem
$$(-\Delta)^s_\Omega u=f\quad{\rm in}\ \, \Omega$$
and show that it is a classical solution when  $f$ is regular. \smallskip

Our  first aim of this paper is to show the existence of solution to  Poisson problem
\begin{equation}\label{eq 4.1-P}
\arraycolsep=1pt\left\{
\begin{array}{lll}
 \displaystyle  (-\Delta)^s_{B_1}   u +u= f\quad   {\rm in}\ \,  B_1,\\[2mm]
\phantom{ ---  }
 \displaystyle   u\in C_0(B_1),
 \end{array}\right.
\end{equation}
when $s\in(0,\frac12]$ and  $f$ satisfies some extra condition.

   \begin{theorem}\label{teo 2}
Assume that $s\in(0,\frac12]$  and
\begin{itemize}
\item[$(\cH_0)$\ ]    $F\in C^\theta_{loc}(\bar B_1 )$ with $\theta\in(0,1)$  is   a non-constant nonnegative function, which is  radially symmetric and  decreasing with respect to $|x|$.
\end{itemize}

Then there exists  $F_0\in \big[\inf_{x\in B_1} F(x),\ \frac1{|B_1|}\int_{B_1} F dx\big)$ such that for $f=F-F_0$,
problem  (\ref{eq 4.1-P}) has a unique classical positive solution $u_{  f}$, which is radially symmetric and decreasing with
respect to $|x|$.

\smallskip

Furthermore, it holds that
$$
\int_{B_1} u_{f}dx=\int_{B_1} f   dx.
$$
  \end{theorem}

\begin{remark}
Note that
$$\int_{B_1} f   dx= \int_{B_1} F dx-F_0 |B_1|>0.$$
We emphasize  that  $f$ can't be a positive constant when it satisfies assumption $(\cH_0)$.
 In fact,  note that if $f=1$  Poisson problem
 $$
    \left\{ \arraycolsep=1pt
    \begin{array}{lll}
    ( - \Delta )_{B_1}^s u +u=1 \quad
     {\rm in} \ \, B_1,\\[2mm]
    \qquad\quad   u\in H^s_0(B_1)
    \end{array}
    \right.
$$
has a unique classical solution $u_1\equiv 1\in H^s_0(B_1)$.  However, it isn't in $C_0(B_1)$.
  In other words, problem (\ref{eq 4.1-P}) has no positive solutions in $C_0(B_1)$ when $f=1$.
\end{remark}

Now we show the existence of positive solutions  to  Schr\"odinger equation (\ref{eq 1.1}).

 \begin{theorem}\label{teo 1}
Assume that  $s\in(0,\frac12]$ and  $$H(d,x,t)=h_1(x)(t-d)^p-d+\epsilon h_2(x),$$
where $p>1$,  $\epsilon>0$ and
 \begin{itemize}
\item[$(\cH_1)$\ ]   functions $h_1, h_2\in C^\theta(\bar B_1)$ with $\theta\in(0,1)$   are  radially symmetric and   decreasing with respect to $|x|$, $h_2$ is non-constant  and
$$ \inf_{x\in B_1}h_1(x),\ \,  \inf_{x\in B_1}h_2(x)>0. $$
\end{itemize}
Denote
$$\epsilon_0=\Big( |B_1| \inf_{x\in B_1} h_2(x) \inf_{x\in B_1}h_1(x)\Big)^{-1} \,  \big(   \|h_1\|_{L^1(B_1)}\big)^{\frac1p}.$$

Then  for any $\epsilon\in(0,\epsilon_0)$, there exists
 $\displaystyle d_\epsilon\in\Big[ \epsilon\inf_{x\in B_1} h_2(x), \ \big(\frac1{ |B_1| \inf_{x\in B_1}h_1(x)}   \|h_1\|_{L^1(B_1)}\big)^{\frac1p}\Big)$
 such that for $h=H(d_\epsilon,\cdot,\cdot)$,
problem
  (\ref{eq 1.1})    admits a   positive solution $u_\epsilon\in C_0(B_1)$ for $\epsilon\in(0,\epsilon_0)$, which is radially symmetric and decreasing with respect to $|x|$.
    \end{theorem}

Note that  for $s\in(\frac12,1)$,  the existence  could  be obtained in $H^s_0(B_1)$ by the variational method,   since the the space $H^s_0(B_1)$ has the boundary trace in \cite{BBC}.  However, it fails  when $s\in (0,\frac12]$.  We emphasize from \cite{CW}  that for $s\in (0,\frac12]$,  Lane-Emden equation
\begin{equation}\label{eq 1.1-LE}
\arraycolsep=1pt\left\{
\begin{array}{lll}
 \displaystyle  (-\Delta)^s_{B_1}   u = u^p\quad & {\rm in}\ \,  B_1,\\[1.5mm]
\phantom{ (-\Delta)^\alpha _{B_1}   }
 \displaystyle   u=0\quad & {\rm on}\   \partial B_1
 \end{array}\right.
\end{equation}
 has no positive solutions.



It is worth noting that  the solutions of (\ref{eq 1.1}) are derived  via passing to the limit of solutions as $r_0\in(0,1)\to 1$ of
  \begin{equation}\label{eq 1.1-MP}
    \left\{ \arraycolsep=1pt
    \begin{array}{lll}
    ( - \Delta )_{B_1}^s u  +u=h(x,u) \quad \
    &{\rm in} \ \, B_{r_0},\\[2mm]
    \qquad\qquad \quad    \ u = 0 \quad \  &{\rm in}\ \  \bar B_1\setminus B_{r_0}.
    \end{array}
    \right.
    \end{equation}
In order to control the boundary behavior in this approximations, we need the special properties of radial  symmetry and the decreasing monotonicity.  Our method for these properties  is  to use  the method of moving planes, which requires some properties of
symmetries and monotonicities for the nonlinearity $h$.

\begin{theorem}\label{teo MP}
Assume that  $s\in(0,1)$,
$$h(x,t)=h_1(x)h_3(t)+h_2(x),$$
where $h_1,\ h_2$ verify
 \begin{itemize}
 \item[$(\cH_2)$\ ]  $h_1,\ h_2:B_1 \to [0,+\infty)$ are  radially symmetric and   decreasing with respect to $|x|$
\end{itemize}
 and $h_3$ satisfies that
 \begin{itemize}
\item[$(\cH_3)$\ ]    $h_3:\R_+\to\R_+$ is nondecreasing,  locally Lipschtiz continuous in $[0,+\infty)$.
\end{itemize}

Let $u\in C_0(B_1)$ be a nonnegative, nonzero solution  of (\ref{eq 1.1-MP}), then $u$ is radially symmetric and strictly decreasing in $r=|x|$ for $r\in(0,r_0)$.
 \end{theorem}


The rest of this paper is organized as follows. In Section 2, we recall the connection of  regional fractional Laplacian and fractional Laplacian, properties of viscosity solutions and  regularity estimates. In Section 3,  we prove Theorem \ref{teo MP} by the method of moving planes.   Section 4 and Section 5 are
devoted to solve the solutions to the related Poisson problem and Schr\"odinger equations, respectively.
  Finally, we annex properties of Green kernel of the fractional Laplacian.




\mysection{Preliminary}

\subsection{Connections of  regional fractional Laplacian and fractional Laplacian}

Note that by the zero extension of the function in $\R^N\setminus \Omega$, we can build the connection between
 regional fractional Laplacian and  fractional Laplacian.

Given $u\in C_0( \Omega  )$, we denote
 \begin{equation}\label{halfe}
  \tilde u (x)=\left\{
      \begin{array}{ll}
      u(x),\quad &   x\in \Omega,\\[2mm]
      0, &    x\notin \Omega.
    \end{array}
    \right.
  \end{equation}
 Then  for $x\in \Omega$, it holds that
  \begin{eqnarray*}
    (-\Delta)^s  \tilde u(x) &=&{\rm p.v.}  \int_{\mathbb{R}^N}\frac{ \tilde u(x)- \tilde u(y)}{|x-y|^{N+2s}}dy\\[2mm]
    &=&{\rm p.v.} \int_{\Omega}\frac{u(x)-u(y)}{|x-y|^{N+2s}}dy+ \int_{\Omega^c}\frac{u(x)}{|x-y|^{N+2s}}dy\\[2mm]
    &=& (-\Delta)^s_\Omega u(x)+ u(x)\varphi_{_\Omega}(x),
  \end{eqnarray*}
where
\begin{equation}\label{ext 1}
\varphi_{_\Omega}(x)=\int_{\Omega^c}\frac{dy}{|x-y|^{N+2s}}.
\end{equation}

\begin{proposition}\label{pr b-ex}
Let $\varphi_{_\Omega}$ be defined by (\ref{ext 1}). \smallskip

$(i)$ If $\Omega$ is  $C^2$, then  $\varphi_{_\Omega}$ is locally Lipschitz continuous.

$(ii)$  If $\Omega=B_1$, then $\varphi_{B_1}$ is radially symmetric, decreasing
and there exists $c_1>0$ such that
 \begin{equation}\label{ext 1-beh}
 \lim_{x\to \partial B_1}\varphi_{B_1}(x)(1-|x|)^{2s}=c_1.
\end{equation}
\end{proposition}
The proof is addressed in the Appendix.

\subsection{ Viscosity solution }

We start with the  definition  of viscosity solutions, inspired by the definition of
viscosity sense for nonlocal problems in \cite{CS}.

\begin{definition}\label{de 2.2}
(i)  We say that a  function  $u\in C( \bar \Omega)$  is a viscosity super-solution (sub-solution)
 of
 \begin{equation}\label{eq 2.1}
 \left\{ \arraycolsep=1pt
\begin{array}{lll}
 \displaystyle  (-\Delta)^{s}_\Omega u=f\quad & {\rm in}\ \,   \Omega,\\[1.5mm]
\phantom{ (-\Delta)^\alpha  }
 \displaystyle   u=0\quad & {\rm on}\    \partial  \Omega,
\end{array}\right.
\end{equation}
if $u\geq 0$ (resp. $u\leq 0$) on $\partial  \Omega$ and
for every point $x_0\in\Omega$ and some  neighborhood $V$ of
$x_0$ with $\bar V\subset \Omega $ and for any $\varphi \in
C^2(\bar V)$ such that $u(x_0)=\varphi(x_0)$ and $x_0$ is the minimum (resp. maximum) point of $u-\varphi$ in $V$,  let
\begin{eqnarray*}
\tilde u  =\left\{ \arraycolsep=1pt
\begin{array}{lll}
\varphi\ \ \ & \rm{in}\ \, &V,\\[1mm]
u \ \ & \rm{in}\ \, &\Omega\setminus V,
\end{array}
\right.
\end{eqnarray*}
we have that
$$(-\Delta)^{s}_\Omega\tilde u(x_0) \geq f(x_0)\quad (resp.  \ (-\Delta)^{s}_\Omega\tilde u(x_0) \le f(x_0)).$$

(ii)  We say that $u$ is a
viscosity solution of (\ref{eq 2.1}) if it is  a viscosity super-solution and also a viscosity sub-solution of (\ref{eq 2.1}).
\end{definition}

\begin{theorem}\label{comparison}
Assume that the functions  $f:\Omega\to\R$, $h:\partial\Omega\to\R$ are
continuous. Let $u$ and $v$ be a viscosity super-solution and
sub-solution  of (\ref{eq 2.1}), respectively. Then
\begin{equation}\label{2.1}
 v  \le u \quad{\rm in}\quad  \Omega.
\end{equation}

\end{theorem}
\noindent{\bf Proof.}  Let us define $w=u-v$, then
\begin{equation}\label{eq 2.2}
 \left\{ \arraycolsep=1pt
\begin{array}{lll}
 \displaystyle  (-\Delta)^{s}_\Omega w \geq 0\quad\, & {\rm in}\ \,  \Omega,\\[1.5mm]
\phantom{ (-\Delta)^s  }
 \displaystyle  w\geq 0\quad\, & {\rm on}\   \partial  \Omega.
\end{array}\right.
\end{equation}
If (\ref{2.1}) fails, then there exists $x_0\in\Omega$ such that
$$w(x_0)=u(x_0)-v(x_0)=\min_{x\in\Omega}w(x)<0,$$
  then in the viscosity sense,
\begin{equation}\label{y 2.1}
 (-\Delta)^{s}_\Omega w(x_0)\geq 0.
\end{equation}
 Since $w$ is a viscosity super solution,  $x_0$ is the minimum point in $\Omega$ and $w\ge 0$ on $\partial\Omega$, then we can
take a small neighborhood $V_0$ of $x_0$ such that $\tilde w=w(x_0)$ in $V_0$, $\tilde w=w$ in $\Omega\setminus V_0$.
From (\ref{y 2.1}), we have that
$$ (-\Delta)^{s}_\Omega \tilde  w(x_0)\geq 0.$$
But the definition of regional fractional Laplacian implies that
$$(-\Delta)^{s}_\Omega \tilde w(x_0)=   \int_{\Omega\setminus V_0} \frac{w(x_0)-w(y)}{|x_0-y|^{N+2\alpha}}dy<0,$$
which is impossible.  \qquad

\begin{remark}\label{rm 2.1-1}
Let $u$ be a continuous function in $\Omega$ and  $x_0$ be a minimum point of $u$, then
$(-\Delta)^{s}_\Omega u(x_0)\geq 0$ in the viscosity sense, where the equality holds if and only if $u$ is a constant.

\end{remark}

  We recall  the stability property for viscosity solutions in our setting.

\begin{theorem}\cite[Theorem 2.2]{CH}\label{stability}
Assume that the function  $g:\Omega\to\R$ is
continuous.  Let $u_n$, $(n\in\N)$ be a sequence of
 functions in  $C(\Omega)$, uniformly  bounded in
$L^1(\Omega)$, $g_n$ and $g$
be continuous in $\Omega$ such that

$(-\Delta)^s_\Omega u_n \ge g_n\ ({\rm resp.}\ (-\Delta)^s_\Omega
u_n \le g_n)$ in $\Omega$ in viscosity sense, $u_n\geq g_n$\ ( resp. \ $u_n\geq g_n)$
on $\partial \Omega$.

$u_n\to u$  locally uniformly in  $\Omega$,

$ u_n\to u$  in  $L^1(\Omega)$,

$ g_n\to g$ locally uniformly in  $\Omega$.\\
    Then $(-\Delta)^s_\Omega  u \ge g\ ({\rm resp.}\ (-\Delta)^s_\Omega u \le g)$ in $\Omega$ in the viscosity sense.
\end{theorem}

\smallskip

Next we have an interior regularity result. For simplicity, we denote by $C^t$ the space $C^{t_0,t-t_0}$ for $t\in(t_0,t_0+1)$, $t_0$ is a positive integer.
\begin{proposition}\label{pr 2.1}
Assume that $s\in(0,1)$, $g\in C_{\rm loc}^{\theta}(\Omega)$ with $\theta>0$, $w\in C_{\rm loc}^{2s+\epsilon}(\cO)\cap L^1(\Omega)$ with $\epsilon>0$
and $2s+\epsilon$  not being an integer, is a solution of
\begin{equation}\label{homo}
(-\Delta)^s_\Omega w=g\quad {\rm in}\ \ \cO .
\end{equation}
Let $\mathcal{O}_1, \mathcal{O}_2$ be  open $C^2$ sets such that
$$\bar{\mathcal{O}_1}\subset\mathcal{O}_2\subset\bar{\mathcal{O}_2}\subset\cO\subset \Omega.$$

Then (i) for any $\gamma\in(0,2s)$ not an integer, there exists $c_2>0$ such that
\begin{equation}\label{2.2}
\norm{w}_{C^\gamma(\mathcal{O}_1)}\le c_2\left(\norm{w}_{L^\infty(\mathcal{O}_2)}+\norm{w}_{L^1(\Omega)}+\norm{g}_{L^\infty(\mathcal{O}_2)} \right);
\end{equation}

(ii) for any $\epsilon'\in(0,\min\{\theta,\epsilon\})$, $2s+\epsilon'$ not an integer, there exists $c_3>0$ such that
\begin{equation}\label{2.3}
\norm{w}_{C^{2s+\epsilon'}(\mathcal{O}_1)}\le c_3\left(\norm{w}_{L^\infty(\mathcal{O}_2)}+\norm{w}_{L^1(\Omega)}+\norm{g}_{C^{\epsilon'}(\mathcal{O}_2)}\right).
\end{equation}

\end{proposition}
{\bf Proof.} The proof is similar to \cite[Proposition 2.1]{CH} for $s\in(\frac12,1)$.  For the reader's convenience, we give the details.  Let $\tilde w=w$ in $\Omega$, $\tilde w=0$ in $\R^N\setminus\bar \Omega$,  we have that
\begin{eqnarray*}
 (-\Delta)^s \tilde w(x)  =  (-\Delta)^s_\Omega w(x)+ w(x)\varphi_{_\Omega}(x),\quad \forall x\in\cO,
\end{eqnarray*}
where $\varphi_{_\Omega}$ is defined as (\ref{ext 1}).
Note that  $\varphi_{_\Omega}\in C^{0,1}_{\rm loc}(\Omega)$.
Combining with (\ref{homo}), we have that
\begin{eqnarray*}
 (-\Delta)^s \tilde w(x)  = g(x)+  w(x)\varphi_{_\Omega}(x),\quad \forall x\in\cO.
\end{eqnarray*}
By \cite[Lemma 3.1]{CV3}, for any $\gamma\in(0,2s)$, we have that
\begin{eqnarray*}
\norm{w}_{C^\gamma(\mathcal{O}_1)}&\le& c_4\left(\norm{w}_{L^\infty(\mathcal{O}_2)}+\norm{w}_{L^1(\Omega)}+
\norm{g+w\varphi_{_\Omega}}_{L^\infty(\mathcal{O}_2)} \right) \\[2mm]
   &\le& c_5\left(\norm{w}_{L^\infty(\mathcal{O}_2)}+\norm{w}_{L^1(\Omega)}+\norm{g}_{L^\infty(\mathcal{O}_2)} \right)
\end{eqnarray*}
 and by  \cite[Lemma 2.10]{RS}, for any $\epsilon'\in(0,\min\{\theta,\epsilon\})$, we have that
 \begin{eqnarray*}
\norm{w}_{C^{2s+\epsilon'}(\mathcal{O}_1)}&\le& c_6\left(\norm{w}_{C^{\epsilon'}(\mathcal{O}_2)}+\norm{g+w\phi}_{C^{\epsilon'}(\mathcal{O}_2)} \right) \\[2mm]
   &\le& c_{7}\left(\norm{w}_{L^\infty(\mathcal{O}_2)}+
   \norm{w}_{L^1(\Omega)}+\norm{g}_{C^{\epsilon'}(\mathcal{O}_2)} \right).
\end{eqnarray*}
We complete the proof.  \hfill$\Box$

  \setcounter{equation}{0}
  \section{Radial symmetry and decreasing monotonicity}

\subsection{ Maximum Principle for small domain}
The essential tool  for the moving planes in a ball is the Maximum Principle for small domains:
\begin{proposition}\label{prop 2-1}
Let $\cO $ be an open  set in ${B_1}$ such that  $|\cO |\leq 2^{-N} |B_{1}|$.
 Suppose that $\phi:\cO \to\R$  is in $L^\infty(\cO )$ satisfying
 \begin{equation}\label{condition for the maximum principle}
\|\phi\|_{L^\infty(\cO )}<+\infty,
\end{equation}
 $w\in L^1(B_1)\cap C(\bar \cO)$  is a solution of
\begin{equation}\label{eq a1}
\left\{ \arraycolsep=1pt
\begin{array}{lll}
  ( - \Delta )_{B_1}^s w \ge  \phi w
  \quad&{\rm in} \
  \cO,\\[2mm]
\qquad \quad \,  w \ge 0\quad \quad \ \ &{\rm in} \ {B_1}\setminus\cO.
\end{array}
\right.
\end{equation}
Then there exists $\delta>0$ such that for $|\cO|<\delta$,
$w\ge0$  in $\cO$.
\end{proposition}

In order to prove Proposition \ref{prop 2-1}, we need the following estimate.

\begin{lemma}\label{lemma 1}
Let $\cO\subset B_1$ be an open set such that  $|\cO |\leq 2^{-N} |B_{1}|$. Suppose that  $g:\cO\to\R$  is in $L^\infty(\cO )$,
 $w \in  L^1(B_1)\cap C(\bar \cO)$ is a solution of
\begin{equation}\label{abp 1}
\left\{ \arraycolsep=1pt
\begin{array}{lll}
( - \Delta )_{B_1}^s w \ge g  \quad&{\rm in} \ \
  \cO,\\[2mm]
\qquad\quad\   w \ge 0\quad& {\rm in} \  \,  B_1 \setminus \cO.
\end{array}
\right.
\end{equation}
 Then  there exists $c_{8}>0$ such that
 \begin{equation}\label{abp}
 -\inf_{\cO}w\le c_{8}\|g\|_{L^\infty(\cO)}|\cO|^{\frac{2s}{N}} .
 \end{equation}

 \end{lemma}

\noindent{\bf Proof.}
The result is obvious if $\inf_{\cO }w\ge0$. Now we assume that $\inf_{\cO}w<0$,
then there exists $x_0\in\cO$ such that
$$w(x_0)=\inf_{x\in\cO}w(x)<0.$$
Combining with (\ref{abp 1}), we have that
\begin{equation}\label{abp471}
-\|g\|_{L^\infty(\cO)}\le g(x_0)
 \le ( - \Delta )_{B_1}^s w(x_0).
 \end{equation}
 By the definition of $( - \Delta )_{B_1}^s$, we have that
 \begin{align*}
  ( - \Delta )_{B_1}^s w(x_0)
&= {\rm p.v.}\int_{B_1}{\frac{w(x_0)-w(y)}{{|{x_0}-y|}^{N+2s}}}dy
\\[2mm]&= {\rm p.v.}\int_{\cO}{\frac{w(x_0)-w(y)}{{|{x_0}-y|}^{N+2s}}}dy
+\int_{B_1\setminus\cO}{\frac{w(x_0)-w(y)}{{|{x_0}-y|}^{N+2s}}}dy
\\[2mm]& \le  \int_{B_1\setminus\cO}{\frac{w(x_0)}{{|{x_0}-y|}^{N+2s}}}dy.
\end{align*}
Let
$$r=\big(|\omega_{_N}|^{-1}|\cO|\big)^{\frac1N}\leq \frac12,$$
 by the   fact that $|\cO |\leq 2^{-N} |B_{1}|$,
it holds that  $|\cO|=|B_r(x_0)|$. We let the vertical plane with respect to $x_0$
$$\cP(x_0)=\big\{z\in\R^N\ |\  z\cdot x_0=0\big\}.$$

Thanks to the decreasing monotonicity of the kernel $\frac{1}{r^{N+2s}}$, we obtain that
\begin{align*}
 \int_{B_1\setminus{B_r(x_0)}}{\frac{1}{|x_0-y|^{N+2s}}dy} \leq \int_{B_1\setminus{\cO }}{\frac{1}{|x_0-y|^{N+2s}}dy},
\end{align*}
since $|B_r(x_0)\setminus \cO|=|\cO \setminus B_r(x_0)|$.
Thus, we derive that
 \begin{align}\label{es-min-1}
 ( - \Delta )_{B_1}^s w(x_0)
 \le  w(0) \int_{B_1\setminus{\cO}}{\frac{1}{|x_0-y|^{N+2s}}dy}
 \le    w(0) \int_{B_1\setminus{B_r(x_0)}}{\frac{1}{|x_0-y|^{N+2s}}dy}.
  \end{align}
 Observe that  for $x_0=0$,  we have that
\begin{align*}
 \int_{B_1\setminus{B_r(x_0)}}{\frac{1}{|x_0-y|^{N+2s}}dy}=\frac1{2s} \omega_{_N} (r^{-2s}-1)\geq \frac1{2s} \omega_{_N} r^{-2s}.
\end{align*}
When $x_0\in B_1\setminus\{0\}$,  let $r_0=\sqrt{1+|x_0|^2}\in(1,\sqrt{2})$,
$$\cP(x_0)= B_1\cap \partial B_{r_0}(x_0)$$
and the cone
   $$\cC(x_0)=\big\{tx_0+(1-t)z:\, \forall\, z\in \cP(x_0),\ \forall\,  t\in(0,1)\big\}. $$
Then
 $\cC(x_0)\subset B_1$,  $|\cC(x_0)|>14 |B_{r_0}(x_0)|$
 and
\begin{align*}
 \int_{B_1\setminus{B_r(x_0)}}{\frac{1}{|x_0-y|^{N+2s}}dy}&> \int_{\cC(x_0)\setminus B_r(x_0)} \frac{1}{|x_0-y|^{N+2s}}dy
 \\[2mm]& \geq\frac14 \int_{B_{r_0}(x_0)  \setminus B_r(x_0)} \frac{1}{|x_0-y|^{N+2s}}dy
  \\[2mm]&=\frac1{8s} \omega_{_N} (r^{-2s}-r_0^{-2s})
 \geq \frac1{8s} \omega_{_N} r^{-2s}.
\end{align*}

 As a consequence,   we derive that
 \begin{align*}
 ( - \Delta )_{B_1}^s w(x_0)
 &\le  w(0) \int_{B_1\setminus{\cO}}{\frac{1}{|x_0-y|^{N+2s}}dy}
  \\[2mm] &\le    w(0) \int_{B_1\setminus{B_r(x_0)}}{\frac{1}{|x_0-y|^{N+2s}}dy}
  \le   \frac1{c_{9}} |\cO|^{-\frac{2s}{N}} w(0),
  \end{align*}
  where $c_{9}=8s|\partial B_1|^{\frac{2s}{N}-1}$.

Finally, together with (\ref{abp471}), we have  that
$$ -\|g\|_{L^\infty(\Omega)}\le ( - \Delta )_{B_1}^s w(x_0)\le  \frac1{c_{9}} w(x_0) |\cO|^{-\frac{2s}{N}}, $$
which implies that
$$w(x_0)\ge -c_{9} \|g\|_{L^\infty(\Omega)}|\cO|^{ \frac{2s}{N}} ,$$
that is,
$$-\inf_{\Omega}w\le c_{9} \|g\|_{L^\infty(\cO)} |\cO|^{\frac{2s}{N}} .$$
We complete the proof. \hfill$\Box$\medskip

 \noindent{\bf Proof of Proposition \ref{prop 2-1}.}  Let us define $\cO^-=\{x\in\cO \ | \ w(x)<0\}$,
 then we observe that
 \begin{equation}\label{eq a1}
 \left\{ \arraycolsep=1pt
 \begin{array}{lll}
 ( - \Delta )_{B_1}^s w(x)\ge \phi(x)w(x),\quad \
&x\in \cO^-,\\[2mm]
\qquad \quad\  w(x)\ge 0,\quad \quad \ \ &x\in {B_1}\setminus{\cO^-}.
\end{array}
\right.
\end{equation}
Using Lemma \ref{lemma 1} with $g =\phi w $, we have that
 \begin{eqnarray*}
\|w\|_{L^\infty(\cO^-)}= -\inf_{\cO^-}w\le c_{9}\|\phi\|_{L^\infty(\cO)}\|w\|_{L^\infty(\cO^-)}|\cO|^{\frac{2s}{N}}.
 \end{eqnarray*}
Then there exists  $\delta>0$  such that for $|\cO|\le \delta$, we have that
$$c_{9} \|\phi\|_{L^\infty(\cO)}  |\cO|^{\frac{2s}{N}}\leq c_{10} \|\phi\|_{L^\infty(\cO)}  \delta^{\frac{2s}{N}}<1,$$
 then $\|w\|_{L^\infty(\cO^-)}=0$, that is, $\cO^-$ is empty.
The proof ends.  \hfill$\Box$ \medskip

\subsection{Moving planes}
  \noindent{\bf Proof of  Theorem \ref{teo MP}. }
Given $\lambda\in (0,r_0)$,  let us define
$$\Sigma_\lambda=\big\{x=(x^1,x')\in B_1\  |\  x^1>\lambda\big\}, \qquad T_\lambda=\big\{x=(x^1,x')\in B_1\ |\  x^1=\lambda\big\},$$
$$\Sigma=\Sigma_\lambda \cup{(\Sigma_\lambda)_{\lambda}},\qquad  w_\lambda(x)=u_\lambda(x)-u(x)$$
and
$$ u_\lambda(x)=\left\{ \arraycolsep=1pt
\begin{array}{lll}
u(x_\lambda),\ \ \ \ &
x\in \Sigma,\\[2mm]
u(x),& x\in  B_1\setminus\Sigma,
\end{array}
\right. $$
where   $x_\lambda=(2\lambda-x^1,x')$  for
$x=(x^1,x')\in B_1.$ For any subset $A$ of $B_1$, we write
$$A_\lambda=\big\{x_\lambda:\, x\in A\big\}.$$ 

{\it Step 1: We  prove that
 $w_\lambda\ge0$ in $\Sigma_\lambda$ if $\lambda\in (0,r_0)$ is close to $r_0$.}
  Indeed, let
 $$\Sigma_\lambda^-=\{x\in \Sigma_\lambda\ |\ w_\lambda(x)<0\}$$
  and
$$ w_\lambda^-(x)=\left\{ \arraycolsep=1pt
\begin{array}{lll}
w_\lambda(x),\ \ \ \ &
x\in \Sigma_\lambda^-,\\[2mm]
0,& x\in  B_1\setminus \Sigma_\lambda^-,
\end{array}
\right.\qquad
w_\lambda^+(x)=\left\{ \arraycolsep=1pt
\begin{array}{lll}
0,\ \ \ \ &
x\in \Sigma_\lambda^-,\\[2mm]
w_\lambda(x),\ \ \ \ & x\in  B_1\setminus  \Sigma_\lambda ^-.
\end{array}
\right. $$
 By the linearity of the regional fractional Laplacian,  we have that  for all $0<\lambda<1$,
$$
 ( - \Delta )_{B_1}^s w_\lambda^+  (x)\le0,\ \ \ \ \forall\ x\in \Sigma_\lambda^-.
$$
In fact, for $x \in  \Sigma_\lambda^-,$  $w_\lambda^+ (x)=0$ and
\begin{eqnarray*}
 ( - \Delta )_{B_1}^s w_\lambda^+ (x) 
 &=&-\int_{{B_1}\backslash  \Sigma_\lambda ^-}  \frac{w_\lambda (y)}{|x-y|^{N+2s}}dy
 \\[2mm]&=&-\int_{{B_1}\backslash {\Sigma}}{\frac{{w_\lambda}(y)}{|x-y|^{N+2s}}}dy-\int_{({\Sigma_\lambda}\setminus \Sigma_\lambda^-) \cup({\Sigma_\lambda}\backslash \Sigma_\lambda^-)^{\lambda}}{\frac{{w_\lambda}(y)}{|x-y|^{N+2s}}}dy
 \\[2mm]&&-\int_{(\Sigma_\lambda^-)^{\lambda}} {\frac{w_\lambda (y)}{|x-y|^{N+2s}}}dy
 \\[2mm]&=:&-I_1-I_2-I_3.
\end{eqnarray*}
Note that
$${I_1}=\int_{{B_1}\setminus {\Sigma}}({u_\lambda}(y)-u(y)){\frac{1}{|x-y|^{N+2s}}}dy=0.$$

Since ${w_\lambda}(y^{\lambda})=-{w_\lambda}(y)$ for any $y \in B_1$, then
\begin{eqnarray*}
 I_2 &=&\int_{({\Sigma_\lambda}\setminus \Sigma_\lambda ^-) \cup({\Sigma_\lambda}\setminus  \Sigma_\lambda ^-)^{\lambda}}{w_\lambda}(y){\frac{1}{|x-y|^{N+2s}}}dy
 \\[2mm]&=&\int_{{\Sigma_\lambda}\setminus  \Sigma_\lambda ^-}{w_\lambda}(y){\frac{1}{|x-y|^{N+2s}}}dy+\int_{{\Sigma_\lambda}\setminus  \Sigma_\lambda ^-}{w_\lambda}(y^{\lambda}){\frac{1}{|x-y^{\lambda}|^{N+2s}}}dy
 \\[2mm]&=&\int_{{\Sigma_\lambda}\setminus  \Sigma_\lambda ^-}{w_\lambda}(y)|Big({\frac{1}{|x-y|^{N+2s}}}-{\frac{1}{|x-{y}^{\lambda}|^{N+2s}}}\Big)dy.
\end{eqnarray*}
For $x \in \Sigma_\lambda^-$ and $y \in {\Sigma_\lambda}\setminus \Sigma_\lambda^-,$ we have
$x-y=(x^1-y^1, x^{'}-y^{'})$, $x-y^{\lambda}=(x^1+y^1-2{\lambda}, x^{'}-y^{'})$,
$|x^1+y^1-2{\lambda}|>|x^1-y^1|$, then
$$\frac{1}{|x-y|^{N+2s}} \ge \frac{1}{|x-{y}^{\lambda}|^{N+2s}}.$$
Combing with $w_{\lambda} \ge0$ in ${\Sigma_\lambda}\setminus  \Sigma_\lambda ^-$, we have that
$${I_2} \ge0.$$
 Since $w_{\lambda}(y)<0$ for $y \in \Sigma_\lambda^-$ and
 ${w_\lambda}(y^{\lambda})=-{w_\lambda}(y)$ for any $y \in B_1$, we have that
\begin{eqnarray*}
 I_3  = \int_{(\Sigma_\lambda^-)^{\lambda}}{w_\lambda}(y){\frac{1}{|x-y|^{N+2s}}}dy
 &=&\int_{ \Sigma_\lambda ^-} {\frac{w_\lambda(y^{\lambda})}{|x-y^{\lambda}|^{N+2s}}}dy
 \\[2mm]&=&-\int_{ \Sigma_\lambda ^-}{\frac{{w_\lambda}(y)}{|x-y^{\lambda}|^{N+2s}}}dy
 \ge  0.
\end{eqnarray*}
Hence, we obtain that for $\lambda\in(0,1)$,
$$( - \Delta )_{B_1}^s w_\lambda^+(x) \le0,\qquad  \forall \, x \in \Sigma_\lambda^-,$$
that is,
$$( - \Delta )_{B_1}^s [w_\lambda-w_\lambda^{-}](x) \le0,\qquad  \forall \, x \in \Sigma_\lambda^-.$$
Then for $x \in \Sigma_\lambda^-$, it holds that
$$( - \Delta )_{B_1}^s w_\lambda (x)\le( - \Delta )_{B_1}^s [w_\lambda^{-}](x)$$
 and
\begin{eqnarray*}
 ( - \Delta )_{B_1}^s  w_\lambda^- (x) & \ge &  ( - \Delta )_{B_1}^s  w_\lambda (x)
\\[2mm]&=& ( - \Delta )_{B_1}^s  u_\lambda  (x)-( - \Delta )_{B_1}^s  u  (x)
\\[2mm]&=& h(x_\lambda,u_\lambda(x))-h(x,u(x))+u(x)-u_\lambda(x)
\\[2mm]&=& \big(h_1(x_\lambda) -h_1(x)\big)h_3(u_\lambda(x))
+h_1(x)\big( h_3(u_\lambda(x))-h_3(u(x))\big)
\\[2mm]&&+ h_2(x_\lambda(x))-h_2(x)+u(x)-u_\lambda(x)
\\[2mm]&\geq & \big(h_1(x)\psi(x)+1\big) \big( u_\lambda(x) - u(x)\big),
\end{eqnarray*}
where the last inequality holds by the assumption $(\cH_2)$ and $(\cH_3)$,
 $$\psi(x)=\frac{h_3(u_\lambda)-h_3(u)}{u_\lambda(x)-u(x)},$$
  which is bounded  for $x\in\Sigma_\lambda^-$.

 Choosing $\lambda\in (0,r_0)$ close enough to $r_0$, then
$|\Sigma_\lambda^-|$ is small enough, by
 $w_\lambda^-=0 \ \ {\rm in}\ \, (\Sigma_\lambda^-)^c$,
  it follows by Proposition \ref{prop 2-1} that
$$w_\lambda=w_\lambda^-\geq0  \ \ {\rm in}\ \,  \Sigma_\lambda^-.$$
 Then $\Sigma_\lambda^-$ is empty, that is,
$$w_\lambda\geq0  \ \ {\rm in}\ \,  \Sigma_\lambda.$$
\smallskip

{\it Step 2:  We claim that for $0<\lambda<1$, if
 $w_\lambda\ge0$ and $w_\lambda\not\equiv0$
 in $\Sigma_\lambda$, then  $w_\lambda>0$
 in $\Sigma_\lambda$.}

If this is not true, then there exists $x_0\in \Sigma_\lambda$ such that
$w_\lambda(x_0)=0$ and then $u_\lambda(x_0)=u(x_0)$ and
\begin{eqnarray}
 ( - \Delta )_{B_1}^s w_\lambda (x_0) &=&( - \Delta )_{B_1}^s   u_\lambda  (x_0)-( - \Delta )_{B_1}^s  u (x_0) \nonumber
 \\[2mm]&=&h((x_0)_\lambda,u_\lambda(x_0))-h(x_0,u(x_0))\nonumber
\\[2mm]&\geq & \big(h_1((x_0)_\lambda -h_1(x_0)\big)h_3(u_\lambda(x_0))\nonumber
\\[2mm]&\geq& 0. \label{e 12-11}
\end{eqnarray}
However,  since $x_0$ is the minimal of $w_\lambda$ and by the definition of the regional
fractional Laplacian
\begin{eqnarray*}
 ( - \Delta )_{B_1}^s [w_\lambda] (x_0) &=&-\int_{B_1}{\frac{{w_\lambda}(y)}{|{x_0}-y|^{N+2s}}}dy
 \\[2mm]&=&-\int_{ \Sigma_\lambda ^-}  {\frac{w_{\lambda} ^-(y)}{|{x_0}-y|^{N+2s}}}dy-\int_{{B_1}\setminus   \Sigma_\lambda ^-} {\frac{w_\lambda^+(y)}{|{x_0}-y|^{N+2s}}}dy
 \\[2mm]&\leq &-\int_{{B_1}\setminus   \Sigma_\lambda ^-} {\frac{w_\lambda^+(y)}{|{x_0}-y|^{N+2s}}}dy
< 0.
\end{eqnarray*}
By the fact that  $ w_ \lambda ^+  \ge0$ and $ w_ \lambda ^+  \not= 0$ in $ \Sigma_\lambda^-$.
Then we obtain a contradiction from (\ref{e 12-11}).
Thus, $w_{\lambda}>0$ in $\Sigma_\lambda$ if ${\lambda} \in(0, 1)$ is close to $r_0$.\smallskip

{\it  Step 3: We show  $\lambda_0=0$, where
$$\lambda_0=\inf\{\lambda\in(0,r_0)\ |\  w_\lambda>0\ \ \rm{in}\ \ \Sigma_\lambda\}.$$}

If it is not true, i.e. $\lambda_0>0$, by the definition of $\lambda_0$, we
have that $w_{\lambda_0}\ge0$ in $\Sigma_{\lambda_0}$ and
$w_{\lambda_0}\not\equiv0$ in $\Sigma_{\lambda_0}$. By {\it Step 2}, we have $w_{\lambda_0}>0$ in $\Sigma_{\lambda_0}$.\smallskip

{\it Claim 1.} If $w_\lambda>0$ in $\Sigma_\lambda$ for
$\lambda\in(0,1)$, then there exists $\epsilon\in(0,\lambda)$ such
that $w_{\lambda_\epsilon}>0$ in $\Sigma_{\lambda_\epsilon}$, where
$\lambda_\epsilon=\lambda-\epsilon$.

Assume that Claim 1 is true, then there exists some
$\epsilon\in(0,\lambda_0)$ such that $w_{\lambda_0-\epsilon}>0$ in
$\Sigma_{\lambda_0-\epsilon}$, which implies that
$$\lambda_0-\epsilon\ge\lambda_0,$$
which is impossible. Then we obtain $\lambda_0=0.$\smallskip

Now we only need to prove Claim 1 to complete  Step 3.\smallskip

 {\it Proof of Claim 1.}
Let $D_\mu=\{x\in\Sigma_\lambda\ | \ dist(x,\partial\Sigma_\lambda)\ge \mu\}$ for $\mu>0$ small. Since
$w_\lambda>0$ in $\Sigma_\lambda$ and $D_\mu$ is compact, then there
exists $\mu_0>0$ such that $w_\lambda\ge \mu_0$ in $D_\mu$. By
continuity of $w_\lambda(x)$, for $\epsilon>0$ small
enough, we denote $\lambda_\epsilon=\lambda-\epsilon,$ then
$$w_{\lambda_\epsilon}(x)\ge0\ \ \rm{in}\ \ D_\mu.$$ As a
consequence, $$\Sigma_{\lambda_\epsilon}^-\subset
\Sigma_{\lambda_\epsilon}\setminus D_\mu$$ and
$|\Sigma_{\lambda_\epsilon}^-|$ small if $\epsilon$ and $\mu$ small.

By  Step 1, $( - \Delta )_{B_1}^s w_{\lambda_\epsilon}^-(x)\le 0$ in $x\in \Sigma_{\lambda_\epsilon}^-$,
Since $w_{\lambda_\epsilon}^+=0$ in
$(\Sigma_{\lambda_\epsilon}^-)^c$ with
$|\Sigma_{\lambda_\epsilon}^-|$ small for  $\epsilon$ and $\mu$
small, $\varphi(x)=\frac{{u_{\lambda_\epsilon}^p}(x)-u^p(x)}{u_{\lambda_\epsilon}(x)-u(x)}$,
 similar with {\it Step 1}, then we have $w_{\lambda_\epsilon}\ge0$ in
$\Sigma_{\lambda_\epsilon}$. And since $\lambda_\epsilon>0$,
$w_{\lambda_\epsilon}\not\equiv0$ in $\Sigma_{\lambda_\epsilon}$,
we have that $w_{\lambda_\epsilon}>0$
 in $\Sigma_{\lambda_\epsilon}$. Thus, Claim 1 is true. \smallskip

 We conclude from the fact of $\lambda_0=0$ that
  $$u(-x^1,x')\ge u(x^1,x')\qquad \ {\rm for}\ \ x^1\ge0.$$
Using the same way, do moving plane from left side to 0, we have
$$u(-x^1,x')\le u(x^1,x')\qquad \ {\rm for}\ \ x^1\ge0.$$
Then $$u(-x^1,x')= u(x^1,x')\qquad \ {\rm for}\ \ x^1\ge0.$$

{\it  Step 4:   we prove $u(x)$ is strictly decreasing in the $x_1$ direction for $x=(x_1,x')\in B_{r_0},$  $x_1>0$.}
By contradiction, if  there exists
$(x^1,x'), (\tilde{x}^1,x')\in \Omega$, $0<x^1<\tilde{x}^1$ such that
\begin{equation}\label{eq 61}
u(x^1,x')\le u(\tilde{x}^1,x').
\end{equation}

Let $\lambda=\frac{x^1+\tilde{x}^1}{2}$ and by arguments above, we have
$$w_\lambda(x)>0\quad \ {\rm for}\ \ x\in\Sigma_{\lambda}.$$
Since $(\tilde{x}^1,x')\in  \Sigma_{\lambda}$, then
\begin{eqnarray*}
0<w_\lambda(\tilde{x}^1,x')
 = u_\lambda(\tilde{x}^1,x')-u(\tilde{x}^1,x')
 = u(x^1,x')-u(\tilde{x}^1,x'),
\end{eqnarray*}
i.e.
$$u((x^1,x'))>u((\widetilde{x}^1,x')),$$
which is impossible with (\ref{eq 61}). Hence, $u(x)$ is strictly decreasing in the
$x_1$ direction for $x=(x^1,x')\in\Omega$ and  $x^1>0$.\hfill$\Box$\medskip


\setcounter{equation}{0}

  \section{Poisson problems}

  In order to prove Theorem \ref{teo 2}, we need the following existence results.

  \begin{lemma} \label{lm 4.1-p}
  Let $s\in(0,1)$, $r \in (0,1)$,  $F: B_1\to [0,+\infty)$ be H\"older continuous, then
  \begin{equation}\label{eq 4.1-P-1}
\arraycolsep=1pt\left\{
\begin{array}{lll}
 \displaystyle  (-\Delta)^s_{B_1}   u +u= F\quad & {\rm in}\ \,  B_{r },\\[2mm]
\phantom{ (-\Delta)^\alpha +u \,  }
 \displaystyle   u=0\quad & {\rm in}\   \bar  B_1\setminus B_{r }
 \end{array}\right.
\end{equation}
has a unique positive solution $u_{r , F}\in C_0(B_1)$.

 Moreover,  (i)  if $F$ is radially symmetric function decreasing with respect to $|x|$,
 then  $u_{r, F}$ is radially symmetric and decreasing with
respect to $|x|$;

 (ii)  $r \to u_{r, F}$ is non-decreasing, i.e.
$$u_{r_1, F}\leq u_{r_2, F}\quad {\rm if}\ \ 0<r_1<r_2<1.$$
  \end{lemma}
  \noindent {\bf Proof.} Let
  $H^s_0(B_{r})$ be the closure of  $C^\infty_0(B_{r})$, with zero value in $\R^N\setminus B_{r}$, under the norm that
  $$\|u\|_{s,r}=\sqrt{\int_{B_1\times B_1} \frac{(u(x)-u(y))^2}{|x-y|^{N+2s}}dxdy+\int_{B_1}u^2 dx},$$
 which is a Hilbert space with the inner product
   $$\langle u,v\rangle _{s,r}= \int_{B_1\times B_1} \frac{(u(x)-u(y))(v(x)-v(y)) }{|x-y|^{N+2s}}dxdy+\int_{B_1}uv dx .$$
 Note that $H^s_0(B_{r})\subset  H^s_0(B_1)$ and from \cite[Corollary 7.2]{EGE},   the embedding $H^s_0(B_1) \hookrightarrow L^q(B_1) $ is
compact for $q\in[2, \frac{2N}{N-2s})$. Since $F$ is H\"older continuous,  it follows by the standard argument of variational methods to find the critical point of
$$\cJ_{s,r}: H^s_0(B_r)\to \R,\qquad  \cJ_{s}(u)=\frac12\|u \|_{s,r}^2-\int_{B_1} F u   dx,  $$
which has a unique critical point  $u_{r, F}\in H^s_0(B_{r})$.  Note that the critical point is the weak solution of (\ref{eq 4.1-P-1})   in the sense that
\begin{equation}\label{e 4.0}
\langle u_{r, F},\xi\rangle _{s,r}=\int_{B_{r}} F \xi dx,\quad\forall\, \xi\in H^s_0(B_{r})
\end{equation}
and taking $\xi=u_{r, F}$, the H\"older inequality and fractional Sobolev embedding  \cite[Theorem 6.7, Remark 6.8]{EGE} implies that
\begin{eqnarray*}
 \|u_{r, F}\|_{s,r}^2&=&\int_{B_{r}} F u_{r, F}  dx
\\[2mm]&\leq & \|u_{r, F} \|_{L^{2^*_s}(B_{r})}\|F\|_{L^{p^*}(B_{r})}
 \leq  c_{11} \|u_{r, F}\|_{s,r}\|F\|_{L^{p^*}(B_1)},
\end{eqnarray*}
that is,
\begin{equation}\label{e 4.1}
 \|u_{r, F}\|_{s,r}  \leq c_{11}  \|F\|_{L^{p^*}(B_1)},
\end{equation}
where   $2^*_s=\frac{2N}{N-2s}$, $p^*=\frac{2N}{N+2s}$ and $ c_{11}>0$ is independent of $r$.

Let
  $\tilde H^s_0(B_{r}) $ be the closure of  $C^\infty_0(B_{r})$, with zero value in $\R^N\setminus B_{r}$, under the norm that
  $$\|\!|u|\!\|_{s,r}=\sqrt{\int_{\R^N \times \R^N} \frac{(u(x)-u(y))^2}{|x-y|^{N+2s}}dxdy+\int_{B_{r}}u^2 dx}$$
  with the inner product
   $$\langle u,v\rangle _{s,r}= \int_{\R^N\times \R^N} \frac{(u(x)-u(y))(v(x)-v(y)) }{|x-y|^{N+2s}}dxdy+\int_{B_{r}}uv dx .$$

Direct computation shows that
$$\|\!|u|\!\|_{s,r}^2=\|u\|_{s,r}^2+ 2\int_{B_1} \varphi_{B_1}(x) u^2(x) dx,$$
where
$\varphi_{B_1}$  is defined in (\ref{ext 1}) with $\Omega=B_1$,
i.e.
$$
\varphi_{B_1}(x)=\int_{B_1^c}\frac{dy}{|x-y|^{N+2s}}.
$$
which by Proposition \ref{pr b-ex} that
$\varphi_{B_1}$ is Lipschitz in $\bar B_{r}$.
Thus, $u_{r, F}$ verifies that
 $$\langle u_{r, F},\xi \rangle _{s,r}={{\int_{B_{r}} F \xi dx}}+ 2\int_{B_{1}} \varphi_{B_1} u_{r, F} \xi dx,
 \qquad \forall\,  \xi\in H^s_0(B_{r}),$$
which means that  $u_{r, F}$ is a weak solution of
\begin{equation}\label{eq 4.1-P-f}
\arraycolsep=1pt\left\{
\begin{array}{lll}
 \displaystyle  (-\Delta)^s   u +u= F+ \varphi_{B_1} u\quad & {\rm in}\ \,  B_{r},\\[1.5mm]
\phantom{ (-\Delta)^\alpha +u \,  }
 \displaystyle   u=0\quad & {\rm in}\ \,    B_{r}^c.
 \end{array}\right.
\end{equation}
Note that the solution could be expressed by
$$u_{r, F}=\Phi_{r}\ast \bar F,  $$
where $\bar F=F+ (\varphi_{B_1}-1) u_{r, F}$ and $\Phi_{r}$ is the Green kernel of $B_{r}$ under the zero condition $\Phi_{r}=0$ in $\R^N\setminus B_{r}$.
Since
$(\varphi_{B_1}-1)$ is uniformly bounded in $B_{r}$ and $u_{r, F}\in L^{2^*_s}$,
then $\bar F\in L^{2^*_s}$ and
if follows by Lemma \ref{embedding-a}
 that $u_{r, F}\in L^\infty(B_{r})$ if $2^*_s>\frac{N}{2s}$
 and we are done;
 or $u_{r, F}\in L^{q_1}(B_{r})$ if
 $q_1=\frac{Nq_0}{2s q_0-N}>q_0$ if   $q_0:=2^*_s>\frac{N}{2s}$,
 in this case, $\bar F\in L^{q_1}$.
 Repeat the procedure, we can find  $i_0\geq 1$ such that $u_{r, F}\in L^\infty(B_{r})$
 and
 $$\|u_{r, F}\|_{L^\infty} \leq c_{12}\|F\|_{L^\infty(B_{r})},    $$
 where $c_{12}>0$ depends on $r$.
 From  (\ref{e 4.1}), we have that $u_{r, F}$ has a uniform bound in $H^s_0(B_{r})$
from Proposition \ref{pr 2.1}, we obtain that for any $\cO\subset \bar\cO\subset B_{r}$,
\begin{eqnarray*}
\norm{u_{r, F}}_{C^\gamma(\cO)} &\leq &c_{13}\left(\norm{u_{r, F}}_{L^\infty(B_{r})}+\norm{u_{r, F}}_{L^1(B_{r})}+{{\norm{F}_{L^\infty(B_{r})}}} \right)
\\[2mm]&\leq &c_{14} \norm{F}_{L^\infty(B_1)}
 \end{eqnarray*}
and then by (\ref{2.3}), we have that
 \begin{eqnarray*}
\norm{u_{r,F}}_{C^{2s+\epsilon'}(\cO)}&\le& c_{15}\left(\norm{u_{r, F}}_{L^\infty(B_{r})}+\norm{u_{r, F}}_{L^1(B_{r})}+{{\norm{F}_{C^\theta(B_{r}) }}} \right)
\\[2mm]&\leq &c_{16} \norm{F}_{C^\theta(B_{r}) }.
\end{eqnarray*}

 It  follows by boundary  regularity in \cite{RS}  that $u_{r, F}$ is a classical solution of (\ref{eq 4.1-P-f}),
 then it is the solution of (\ref{eq 4.1-P-1}).  Thus, for $\xi\in C^{2s+\theta}_{loc}(B_{r})\cap C^s_0(B_{r})$,
 it holds that
 \begin{eqnarray*}
\int_{B_{r}}F\xi dx  &=& \int_{B_{r}} \big(\xi (-\Delta)^s_{B_1} u_{r, F} +u_{r, F} \xi\big) dx
 \\[2mm] &=& \int_{B_{r}} \big(\xi (-\Delta)^s  u_{r, F} +u_{r, F} \xi-\varphi_{B_1}(x)u_{r, F} \xi \big) dx
 \\[2mm]&=&\int_{B_{r}} \big( u_{r, F} (-\Delta)^s \xi+u_{r, F} \xi-\varphi_{B_1}(x)u_{r, F} \xi \big) dx
  \\[2mm]&=&\int_{B_{r}} \big(u_{r, F} (-\Delta)^s_{B_1}  \xi+u_{r, F} \xi\big) dx,
 \end{eqnarray*}
 that is,
\begin{equation}\label{e 4.1-dis}
\int_{B_{r}} u_{r, F}\big( (-\Delta)^s_{B_1}  \xi+  \xi\big) dx =\int_{B_{r}}F\xi dx.
\end{equation}
 By the method of moving planes in Theorem \ref{teo MP} with $h_1=0$ and {{$h_2=F$}}, the solution {{$u_{r, F}$ }}is radially symmetric and
 decreasing to $|x|$.

 For $0<r_1<r_2<1$,  we see that the solution $u_{r_2, F}$ is a super solution of (\ref{eq 4.1-P-1})
 with $r=r_1$, then the maximum principle shows that $u_{r_2, F}\geq u_{r_1, F}$ in $B_1$.
   \hfill$\Box$\medskip

  \noindent{\bf Proof of Theorem \ref{teo 2}.}
   It follows by Lemma \ref{lm 4.1-p}, problem (\ref{eq 4.1-P-1}) has a unique classical solution $u_{r,F}$, which is positive, radially symmetric and decreasing with respect to $|x|$ and $r:\to u_{r,F}$ is increasing. From (\ref{e 4.1})
  $$\|u_{r,F}\|_{s, r}\leq c_{17}\|F\|_{L^{p^*}(B_1)}, $$
  where $p^*= \frac{2N}{N+2s}$ and $c_{17}>0$ is independent of $r$.
  Note that $r\to u_{r,F}$ is increasing,  passing to the limit as $r\to1^-$, it yields that
  $$u_{1,F}=\lim_{r\to1^-}u_{r,F}\quad {\rm as}\ \, r\to1^-\quad {\rm weakly\ in\ } H^s_0(B_1)\quad
   {\rm and\ \  a.e.\ in} \ B_1. $$
 By compact embedding, the above convergence holds   strongly  in  $L^{p}(B_1)$  for $p\in [2,2^*_s)$.
  From (\ref{e 4.0}), we have that
  \begin{equation}\label{e 4.0-r=1}
\langle u_{1,F},\xi\rangle _{s,1}=\int_{B_1} F \xi dx,\qquad\forall \, \xi\in H^s_0(B_1)
\end{equation}
and for $\xi\in   C^2_c(B_1)$,
\begin{equation}\label{e 4.1-dis r=1}
\int_{B_1} u_{1,F}\big( (-\Delta)^s_{B_1}  \xi+  \xi\big) dx =\int_{B_1}F\xi dx.
\end{equation}
The function {{$u_{1,F}$}} is  the critical point of the functional
$$\cJ_{s}: H^s_0(B_1)\to \R,\qquad  \cJ_{s}(u)=\frac12\|u \|_{s,1}^2-\int_{B_1} F u   dx,  $$
whose critical point is unique.
Then $ \int_{B_1} u_{1,F}  dx =\int_{B_1}F  dx$ by taking $\xi \equiv 1\in H^s_0(B_1)$ in (\ref{e 4.0-r=1}) for $s\in(0,\frac12]$.

  Moreover,   $u_{1,F}$ inherits the positivity, the symmetry property and decreasing monotonicity of $u_{r,F}$, then $u_{1,F}$ is locally bounded in $B_1\setminus\{0\}$.

To prove  $u_{1,F}\in L^\infty(B_1)$.  Let $t>1$  and
 $$w_t =(u_{1,F}-t)_+\quad {\rm in}\ \, B_1, $$
 which is $H^s_0(B_1)$.
For $\sigma\in(0,\frac12]$, we have  that
  \begin{eqnarray*}
 \int_{B_1}{{ F w_t }}dx &=& \langle {{u_{1,F},}}w_t\rangle _{s}
 \\[1mm] &=&\int\int_{\{x,y\in B_1:|x-y|<\sigma\}}\frac{\big(u_{1,F}(x)-u_{1,F}(y)\big)\big(w_t(x)-w_t(y)\big)}{|x-y|^{N+2s}}dxdy
 \\[1mm]&&+\int\int_{\{x,y\in B_1:|x-y|\geq \sigma\}}\frac{\big(u_{1,F}(x)-{{u_{1,F}(y)}}\big)\big(w_t(x)-w_t(y)\big)}{|x-y|^{N+2s}}dxdy
   \\[1mm]&\geq&\int\int_{\{x,y\in B_1:|x-y|<\sigma\}}\frac{ \big(w_t(x)-w_t(y)\big)^2}{|x-y|^{N+2s}}dxdy
   \\[1mm]&&\,+2c_{18}\Big(\sigma^{-2s}\int_{B_1} u_{1,F} w_t dx
    -\int_{B_1} (\kappa_\sigma\ast u_{1,F}) w_t dx\Big)
  \\[1mm]&:=& \cF_\sigma+2c_{18}(\cE_{1,\sigma}-\cE_{2,\sigma}),
\end{eqnarray*}
where $c_{18}=\frac{\omega_{N}}{2s}$ and $\kappa_\sigma =\chi_{\R^N\setminus B_\sigma} |\cdot|^{-N-2s}$,
the last inequality holds by the fact that
   \begin{eqnarray*}
&&   \big(u_{1,F}(x)-u_{1,F}(y)\big)\big(w_t(x)-w_t(y)\big)
\\[2mm]&=&\big((u_{1,F}(x)-t)-(u_{1,F}(y)-t)\big)\big(w_t(x)-w_t(y)\big)
    \\[2mm]&=&(u_{1,F}(x)-t)w_t(x)+(u_{1,F}(y)-t)w_t(y)-(u_{1,F}(x)-t)w_t(y)-(u_{1,F}(y)-t)w_t(x)
     \\[2mm]&=&w_t^2(x)+w_t^2(y)-2w_t(x)w_t(y)+(u_{1,F}(x)-t)_-w_t(x)+(u_{1,F}(y)-t)_-w_t(y)
      \\[2mm]&\geq&\big(w_t(x)-w_t(y)\big)^2
 \end{eqnarray*}
 for  $t_-=\max\{-t,0\}$ and $x,y\in B_1$.

 Direct computations show that
    \begin{eqnarray*}
    \|\kappa_\sigma\ast u_{1,F}\|_{L^\infty(B_1)}&\leq &\|u_{1,F}\|_{L^{2^*_s}(B_1)}\big(\int_{\R^N\setminus B_\sigma} |y|^{-(N+2s)p^* }ds\big)^{\frac1{p^*}}
      \\[2mm]&\leq&c_{19} \|F\|_{L^{p^*}(B_1)}  \sigma^{-\frac{N+2s}{2}}
       \\[2mm]&\leq&c_{20} \|F\|_{L^{\infty}(B_1)}\sigma^{-\frac{N+2s}{2}},
  \end{eqnarray*}
  which implies that
  $$\cE_{2,\sigma}\leq c_{20} \|F\|_{L^{\infty}(B_1)}\sigma^{-\frac{N+2s}{2}}  \int_{B_1} w_t dx.$$
  Moreover, it holds that
 $$ \int_{B_1}{{ F w_t dx}}\leq \|F\|_{L^\infty(B_1)} \int_{B_1} w_t dx, \qquad   \cE_{1,\sigma}\geq c_{21}\sigma^{-2s} t\int_{B_1} w_t dx.$$
 As a consequence, if $t\geq t_0$ for some $t_0>0$ large enough,  we obtain that
 $$\cF_\sigma\leq \big( \|F\|_{L^{\infty}(B_1)}(c\sigma^{-\frac{N+2s}{2}}+1) -c_{21}\sigma^{-2s}t\big)\int_{B_1} w_t dx\leq0, $$
 where
  $$  \|F\|_{L^{\infty}(B_1)}(c\sigma^{-\frac{N+2s}{2}}+1) -c_{21}\sigma^{-2s}t <0\quad {\rm if}\ \, t\ {\rm is\ large}. $$
 So we have that $w_t=0$ a.e. in $\Omega$, which means $u_{1,F}\leq t_0$.  The $L^\infty$ bound of $u_{1,F}$ is obtained. \smallskip

  Note that for $\cO\subset\bar \cO \subset B_1$, then for $r$ close to 1 such that $\cO\subset\bar \cO \subset B_{r}$,
   \begin{eqnarray*}
{{\norm{u_{r, F}}}}_{C^\gamma(\cO)} &\leq &c_{22}\left(\norm{u_{1,F}}_{L^\infty(B_1)}+\norm{u_{1, F}}_{L^1(B_{r})}+\norm{F}_{L^\infty(B_1)} \right)
\\[2mm]&\leq &c_{23} \big(\norm{F}_{L^\infty(B_1)}+\norm{u_{1,F}}_{L^\infty(B_1}\big)
 \end{eqnarray*}
 and  by (\ref{2.3}), we have that
 \begin{eqnarray*}
\norm{u_{r, F}}_{C^{2s+\epsilon'}(\cO)}&\le& c_{24}\left(\norm{u_{1,F}}_{L^\infty(B_1)}+\norm{w}_{L^1(B_{r})}+\norm{F}_{C^\theta(B_{r}) } \right)
\\[2mm]&\leq &c_{25} \norm{F}_{C^\theta(B_1) }.
\end{eqnarray*}
Then $u_{r,F}\to u_{1,F}$ as $r\to1^-$ locally in  $C^{2s+\epsilon''}(B_1)$.
By the stability results of Theorem \ref{stability}, we obtain that $u_{1,F}\in C^{2s+\epsilon}(B_1)$ and it verifies that
 $$(-\Delta)^s_{B_1} u_{1,F}+u_{1,F}=F\quad{\rm in}\ \, B_1. $$


As proved above,  $u_{1,F}$ is   radially symmetric function decreasing with respect to $|x|$,  then we can denote
$$d_{1,F}=\lim_{|x|\to1^-} u_{1,F}(x)\geq0. $$
Let
$$ u_f(x)=u_{1,F}-d_{1,F},$$
then $u_f\in C_0(B_1)$ is nonnegative   and verifies that
 \begin{equation}\label{eq 4.1-P-f=1}
\arraycolsep=1pt\left\{
\begin{array}{lll}
 \displaystyle  (-\Delta)^s_{B_1}  u_f +u_f = F-d_{1,F} \quad & {\rm in}\ \,  B_1,\\[2mm]
\phantom{ (-\Delta)^\alpha_{B_1} +u_f    }
 \displaystyle  u_f=0\quad & {\rm on}\   \partial B_1.
 \end{array}\right.
 \end{equation}
and
$$0\leq \int_{B_1}u_fdx=\int_{B_1} (F-d_{1,F})dx. $$
If $F=f_0$ is a constant, then $f_0-d_{1,F}$ is a unique solution of
$$
\arraycolsep=1pt\left\{
\begin{array}{lll}
 \displaystyle  (-\Delta)^s_{B_1}   u +u= F-d_{1,F}\quad   {\rm in}\ \,  B_1,\\[2mm]
\phantom{ (-\Delta)^\alpha +u \,  }
 \displaystyle   u\in H^s_0(B_1)
 \end{array}\right.
$$
  and by the uniqueness,  $u_f\equiv f_0-d_{1,F}$ and by the zero boundary, we have that $f_0=d_{1,F}$.

If $F$ is not a constant,  $u_f$ is no longer a constant, $\int_{B_1}u_fdx>0$,   then $\int_{B_1} (F-d_{1,F})dx>0$
which implies  $d_{1,F}<\frac{1}{|B_1|} \int_{B_1} F(x)dx$.

Finally, we claim that $d_{1,F}\geq \inf_{x\in B_1} F(x)$.   In fact,  Letting  $d_0:=\inf_{x\in B_1} F(x)>0$,  then by comparison principle
$$u_{r,d_0}\leq u_{r, F}\quad {\rm in }\ B_1$$
where $u_{r, d_0}$, $u_{r, F}$ are
the solutions  of (\ref{eq 4.1-P-1}) with non-homogeneous term $d_0$  and  $F$ respectively.

By the convergence, we obtain that
$$u_{1,F}\geq u_{1,d_0}\equiv d_0\quad {\rm in}\ B_1, $$
which implies that $d_{1,F}\geq \inf_{x\in B_1} F(x).$  \hfill$\Box$\medskip



 From the the proof of  Theorem \ref{teo 2}, we conclude that
   \begin{corollary}\label{cr p-xxx}
Assume that  $s\in(0,\frac12]$,  $F\in C^\theta(\bar B_1)$ with $\theta\in(0,1)$, is a  nonnegative function,  radially symmetric and  decreasing with respect to $|x|$.

Then    problem
\begin{equation}\label{eq 4.1-P-xx}
\arraycolsep=1pt\left\{
\begin{array}{lll}
 \displaystyle  (-\Delta)^s_{B_1}   u +u= F\quad   {\rm in}\ \,  B_1,\\[2mm]
\phantom{ (- \,  }
 \displaystyle   u\in H^s_0(B_1)
 \end{array}\right.
\end{equation}
has a unique positive solution $u_{1 , F}\in C(\bar B_1)$, which is radially symmetric and decreasing with
respect to $|x|$.

Moreover, (i) the mapping: $F\mapsto u_{1,F}$ is increasing and
$$d_{1,F}\in\Big[\inf_{x\in B_1} F(x),\frac{1}{|B_1|}\int_{B_1} Fdx\Big];$$

(ii) if $F$ is a positive constant, we derive that $u_{1,F}=F$.
  \end{corollary}
{\bf Proof. } For   $F_1\leq F_2$,  $ u_{r,F_1}\leq u_{r,F_2}$ by the previous proof, passing to the limit we get
the mapping: $F\mapsto u_{1,F}$ is increasing.

If $F$ is a constant, then $w:=u_{1,F}-F$ is a solution of
$$
\arraycolsep=1pt\left\{
\begin{array}{lll}
 \displaystyle  (-\Delta)^s_{B_1}   u +u= 0\quad   {\rm in}\ \,  B_1,\\[2mm]
\phantom{ (- \,  }
 \displaystyle   u\in C_0(B_1)
 \end{array}\right.
$$
which only has a zero solution by the maximum principle. Then  $u_{1,F}=F$.  \hfill$\Box$\medskip

  \setcounter{equation}{0}
  \section{Schr\"odinger equation}
Under the assumption of Theorem \ref{teo 1},    Schr\"odinger equation (\ref{eq 1.1}) could be written as
 \begin{equation}\label{eq 4.1-S}
\arraycolsep=1pt\left\{
\begin{array}{lll}
 \displaystyle  (-\Delta)^s_{B_1}   u +u=h_1  u^p+\epsilon h_2  \quad  {\rm in}\ \,  B_1,\\[2mm]
\phantom{ (-=    }
 \displaystyle   u\in H^s_0(B_1),
 \end{array}\right.
\end{equation}
where $p>1$ and $\epsilon>0$. \smallskip

\noindent{\bf Proof of Theorem \ref{teo 1}.} Let $u_{h_2}$ be the unique solution of
 $$ (-\Delta)^s_{B_1}   u +u=  h_2  \quad  {\rm in}\ \,  B_1, \qquad    u\in H^s_0(B_1). $$
 Now we define the iterating sequence
$$v_0:=\epsilon  u_{h_2}>0,$$
and by Corollary \ref{cr p-xxx}, $v_n$ with $n=1,2,\cdots$ is the unique solution of
  \begin{equation}\label{eq 4.1-s-n}
  (-\Delta)^s_{B_1}   u +u= h_1 v_{n-1}^p+ \epsilon h_2 \quad  {\rm in}\ \,  B_1, \qquad    u\in H^s_0(B_1).\end{equation}
and we have that $v_1\geq v_0$.
Assuming that
$$
v_{n-1} \ge  v_{n-2} \quad{\rm in} \quad B_1,
$$
 then
\begin{eqnarray*}
  (-\Delta)^s_{B_1}  (v_n-v_{n-1})  + (v_n-v_{n-1}) =h_1( v_{n-1}^p-v_{n-2}^p) \geq0\quad {\rm in}\ \, B_1
\end{eqnarray*}
and
$v_n-v_{n-1}\in H^s_0(B_1),$ we apply Corollary \ref{cr p-xxx} to obtain that $v_n\geq v_{n-1}$ in $B_1$.

Thus the sequence $\{v_n\}_{n\in\N}$ is  increasing with respect to $n$.

We next build an upper bound for the sequence $\{v_n\}_n$.  For  $t>0$, denote
$$
 w_t= t,
$$
 then
 \begin{eqnarray*}
  (-\Delta)^s_{B_1}  w_t +w_t-h_1 w_t^p
= t -t^p h_1
\geq t  \Big(1-t^{p-1}\|h_1  \|_{L^\infty(B_1)} \Big)
\end{eqnarray*}
and letting
$$L(t)= t   -t^p\|h_1  \|_{L^\infty(B_1)} ,$$
note that $L(\cdot)$ has maximum $\frac{p-1}{p} \big(p\|h_1  \|_{L^\infty(B_1)}\big)^{-\frac{1}{p-1}}$  in
 at $t_p=(p\|h_1 \|_{L^\infty(B_1)})^{-\frac{1}{p-1}}$.

In order to find the upper solution, we take $t=t_p$ and if
\begin{equation}\label{wt}
 \epsilon \|h_2\|_{L^\infty(B_1)} \leq  \frac{p-1}{p} \big(p\|h_1\|_{L^\infty(B_1)} \big)^{-\frac{1}{p-1}}.
\end{equation}
then
\begin{equation}\label{4.2.4}
 (-\Delta)^s_{B_1}  w_{t_p} +w_{t_p}\geq h_1 w_{t_p}^p+\epsilon h_2.
\end{equation}
Note that (\ref{wt}) holds if
$$\epsilon\leq \epsilon_p:=\frac{p-1}{p} \big(p\|h_1 \|_{L^\infty(B_1)} \big)^{-\frac{1}{p-1}}\|h_2\|_{L^\infty(B_1)} ^{-1} $$

Obviously, we have that {{$ w_{t_p } \geq v_0$}}.
Inductively, we obtain
\begin{equation}\label{2.10a}
v_n\le {{ w_{t_p}}}
\end{equation}
for all $n\in\N$. Therefore, the sequence $\{v_n\}_n$ converges. Let $u_{\epsilon}:=\displaystyle \lim_{n\to\infty} v_n$ in $B_1$.  By the regularity results,  $u_{\epsilon}$ is a  solution of (\ref{eq 4.1-S}).\smallskip

We claim that $u_{\epsilon}$ is the minimal solution of (\ref{eq 1.1}), that is, for any nonnegative solution $u$ of (\ref{eq 1.2}), we always have $u_{\epsilon}\leq u$. Indeed,  there holds
\[
  (-\Delta)^s_{B_1}   u  +  u     = h_1  u^p +\epsilon h_2 \ge  (-\Delta)^s_{B_1}   v_0+v_0\quad {\rm in}\ \, B_1,\qquad  u_{\epsilon}=u\quad {\rm on}\ \, \partial B_1
\]
 then $u\geq v_0$, $u_{\epsilon}=u$ on $\partial B_1$ and
 \[
  (-\Delta)^s_{B_1}   u  +  u     = h_1  u^p +\epsilon h_2 \geq h_1 v_0^p +\epsilon h_2  =  (-\Delta)^s_{B_1}   v_1+v_1\quad {\rm in}\ \, B_1,
\]
which  implies that  $u\geq v_1$ in $B_1$.
 We may show inductively that
\[
u\ge v_n
\]
for all $n\in\N$.  The claim follows.

From above argument,   if problem (\ref{eq 4.1-S}) has a nonnegative solution $u_{\epsilon_1}$  for $ \epsilon_1>0$, then (\ref{eq 4.1-S}) admits a minimal solution $u_{\epsilon}$ for all $ \epsilon\in(0, \epsilon_1]$. As a result, the mapping $ \epsilon\mapsto u_{\epsilon}$ is increasing.
So we may define
$$\epsilon^*=\sup\big\{\epsilon>0:\ (\ref{eq 4.1-S})\ {\rm has\ minimal\ solution\ for\ } \epsilon \big\}$$
and we have that
$$\epsilon^*\ge \epsilon_p.$$

Finally,  we  prove that $\epsilon^*<+\infty$. Assume that (\ref{eq 4.1-S}) has a positive solution
    for $\epsilon>0$. Our above proof shows that (\ref{eq 4.1-S})
  has a minimal solution $u_\epsilon$. Let   $u_{h_1}$
  be the solution of
  $$
  \arraycolsep=1pt\left\{
\begin{array}{lll}
 \displaystyle  (-\Delta)^s_{B_1}   u +u= h_1\quad   {\rm in}\ \,  B_1,\\[2mm]
\phantom{ (-\,  }
 \displaystyle  u\in H^s_0(B_1).
 \end{array}\right.
 $$
 If $h_0:=\inf_{x\in B_1} h_1(x)>0$, then $u_{h_1}\geq h_0$ by Corollary \ref{cr p-xxx}.

Letting $u_{h_1}$ as test function, we have that
 \begin{eqnarray*}
 \int_{B_1}   u_{\epsilon}^p h_1 u_{h_1} dx+\epsilon \int_{B_1} h_2 u_{h_1}dx
 &=&\int_{B_1}  \Big( (-\Delta)^s_{B_1}   u_{\epsilon}  +  u_{\epsilon}\Big) u_{h_1} dx
 \\[2mm]  &=& \int_{B_1} u_{\epsilon} \Big(  (-\Delta)^s_{B_1}   u_{h_1} +u_{h_1} \Big) dx
 \\[2mm] &=& \int_{B_1} u_{\epsilon}h_1dx
  \\[2mm] &\leq& \Big(\int_{B_1} u_{\epsilon}^p h_1  u_{h_1} \,dx\Big)^{\frac1p}  \Big(\int_{B_1} h_1 (u_{h_1})^{-\frac1{p-1}} \,dx\Big)^{1-\frac1p}
   \\[2mm] &\leq & c_{26}\Big(\int_{B_1} u_{\epsilon}^p h_1 u_{h_1} \,dx\Big)^{\frac1p},
 \end{eqnarray*}
where
 \begin{eqnarray*}
 c_{26}= \big(\int_{B_1}  h_1(u_{h_1})^{-\frac1{p-1}} \,dx\big)^{1-\frac1p}
  \leq h_0^{1-\frac1p}  \|h_1\|_{L^1(B_1)} ^{1-\frac1p}
 <+\infty.
 \end{eqnarray*}
 Thus, we have that
\begin{equation}\label{2.14-x}
\int_{B_1} u_{\epsilon}^ph_1u_{h_1} \,dx \le  c_{26}^{\frac{p}{p-1}}
\end{equation}
and we have that
\begin{equation}\label{2.14}
\epsilon \le  \frac{ c_{26}^{\frac{p}{p-1}}}{\int_{B_1} h_2 u_{h_1}dx}=\frac{\int_{B_1}  h_1 (u_{h_1})^{-\frac1{p-1}} \,dx}{\int_{B_1} h_2 u_{h_1}dx},
\end{equation}
which means
$$\epsilon^*\leq \frac{\int_{B_1}  h_1 (u_{h_1})^{-\frac1{p-1}} \,dx}{\int_{B_1} h_2 u_{h_1}dx}<+\infty.$$

Finally,   $u_\epsilon$ is  radially symmetric function decreasing with respect to $|x|$,  then we can denote
$$d_{\epsilon}=\lim_{|x|\to1^-} u_{\epsilon}(x)\geq0. $$
By by Corollary \ref{cr p-xxx}, we have that
$$d_{\epsilon}\geq \epsilon\lim_{|x|\to1^-} h_{h_2}\geq \epsilon\inf_{x\in B_1} h_2(x).$$
By (\ref{2.14-x}), we have that
\begin{align*}
h_0^2 |B_1| d_{\epsilon}^p&\leq   d_{\epsilon}^p \int_{B_1} {{h_1u_{h_1} }}dx\\[2mm]
& \leq \int_{B_1} u_{\epsilon}^ph_1u_{h_1} \,dx
  \le  c_{26}^{\frac{p}{p-1}} \leq  h_0   \|h_1\|_{L^1(B_1)},
\end{align*}
 that is
$$d_{\epsilon} \leq \Big(\frac1{h_0|B_1|}   \|h_1\|_{L^1(B_1)}\Big)^{\frac{1}{p}}.$$

Let
$$w_\epsilon (x)=u_\epsilon-d_{\epsilon},$$
then $w_\epsilon\in C_0(B_1)$ is nonnegative   and verifies that
 \begin{equation}\label{eq 4.1-P-fxxxx}
\arraycolsep=1pt\left\{
\begin{array}{lll}
 \displaystyle  (-\Delta)^s_{B_1}  w_\epsilon +w_\epsilon =
  h_1  (w_\epsilon+d_{\epsilon})^p-d_{\epsilon}+\epsilon h_2  \quad & {\rm in}\ \   B_1,\\[2mm]
\phantom{ (-\Delta)^\alpha_{B_1} +u_f    }
 \displaystyle  w_\epsilon=0\quad & {\rm on}\  \    \partial B_1.
 \end{array}\right.
 \end{equation}
and
$$0\leq \int_{B_1}w_\epsilon dx=\int_{B_1} \Big( h_1  (w_\epsilon+d_{\epsilon})^p-d_{\epsilon}+\epsilon h_2\Big)dx. $$
 The proof ends.    \hfill$\Box$

 \setcounter{equation}{0}
  \appendix
 \section{ Some estimates }
 \subsection{Proof of Proposition \ref{pr b-ex} }
(i)  For
 $x_1,x_2\in \Omega $ and  any $z\in \R^N\setminus \Omega$, we have that
 $$|z-x_1|\geq \rho(x_1)+\rho(z), \qquad |z-x_2|\geq \rho(x_2)+\rho(z)$$ and
$$||z-x_1|^{N+2s}-|z-x_2|^{N+2s}|\leq
c_{27}|x_1-x_2|(|z-x_1|^{N+2s-1}+|z-x_2|^{N+2s-1}),
$$
where $\rho(x)={\rm dist}(x, \partial \Omega)$,   $c_{27}>0$ is independent of $x_1$ and $x_2$. Then
\begin{align*}
 &\quad |\varphi_{_\Omega}(x_1)-\varphi_{_\Omega}(x_2)| \\[2mm]
 &\leq \int_{  \Omega^c }  \frac{||z-x_2|^{N+2s}-|z-x_1|^{N+2s}|}{|z-x_1|^{N+2s}|z-x_2|^{N+2s}}dz\\[2mm]
 &\leq c_{27} |x_1-x_2|\left[\int_{ \Omega^c}  \frac{dz}{|z-x_1||z-x_2|^{N+2s}}+\int_{ \Omega^c}  \frac{dz}{|z-x_1|^{N+2s}|z-x_2|}\right].
\end{align*}
By direct computation, we have that
\begin{align*}
 \int_{ \Omega^c}  \frac{1}{|z-x_1||z-x_2|^{N+2s}}dz
  &\le  \int_{\R^N\setminus{B_{\rho(x_1)}(x_1)}} \frac{1}{|z-x_1|^{N+2s+1}}dz
 \\[2mm]&\qquad +\int_{\R^N\setminus{B_{\rho(x_2)}(x_2)}} \frac{1}{|z-x_2|^{N+2s+1}}dz \\[2mm]
 &\le  c_{28}(\rho(x_1)^{-1-2s}+\rho(x_2)^{-1-2s})
\end{align*}
and similar to obtain that
$$\int_{ \Omega^c} \frac{1}{|z-x_1|^{N+2s}|z-x_2|}dz\le  c_{29}(\rho(x_1)^{-1-2s}+\rho(x_2)^{-1-2s}),$$
where $c_{28},c_{29}>0$ are independent of $x_1$ and $ x_2$.
Then
$$|\varphi_{_\Omega}(x_1)-\varphi_{_\Omega}(x_2)|
\le  c_{30}(\rho(x_1)^{-1-2s}+\rho(x_2)^{-1-2s})|x_1-x_2|,$$
where $c_{30}=c_{27} (c_{28}+c_{29})$,
it implies that $\varphi_{_\Omega}$ is locally Lipschitz continuous.\smallskip

 (ii)  Firstly, we claim that $\varphi_{B_1}(x)=\varphi_{B_1}(z)$ if $|x|=|z|$. In fact,
denote  ${\bf A}$  a matrix with $|{\bf A}|=1$ and $z={\bf A}x$, we have that
 \begin{eqnarray*}
  \varphi_{B_1}(z)  = \varphi_{B_1}({\bf A}x)
 &=&\int_{B_1^c}\frac{dy}{|{\bf A}x-y|^{N+2s}}
  \\&=&\int_{B_1^c}\frac{d\tilde{y}}{|x-\tilde{y}|^{N+2s}}
  = \varphi_{B_1}(x),
\end{eqnarray*}
where $\tilde y={\bf A}^{-1} y$.

Now we show the monotonicity.  By the radial symmetry of $\varphi$, we let
 $$\varphi(r)=\varphi_{B_1}(x),  \quad \ r=|x|\in (0,1).$$
Fixed $x_{1}=t_{1}e_{1},$ $x_{2}=t_{2}e_{1}$,
$e_{1}=(1,0,\cdots,0)\in \mathbb{R}^N$,  $0<t_{1}<t_{2}<1$, by direct computation, it yields that
 \begin{align*}
  \varphi(t_1)-\varphi(t_2) &= \int_{B_1^c}(\frac{1}{|t_{1}e_{1}-y|^{N+2s}}-\frac{1}{|t_{2}e_{1}-y|^{N+2s}})dy
  \\[2mm]&= \int_{\cA_{1}\cup\cA_{2}}(\frac{1}{|t_{1}e_{1}-y|^{N+2s}}-\frac{1}{|t_{2}e_{1}-y|^{N+2s}})dy
  \\[2mm]& \quad\  +\int_{\cA_0}(\frac{1}{|t_{1}e_{1}-y|^{N+2s}}-\frac{1}{|t_{2}e_{1}-y|^{N+2s}})dy,
\end{align*}
where $\cA_0= B_1\big((t_1+t_2)e_1\big)\setminus B_1,$
$$\cA_{1}=\left\{(x^1,x')\  | \ (x^1,x')\in \big(-\infty, \frac{t_1+t_2}{2}\big)\times \R^{N-1}\setminus B_1 \right\} $$
and $$ \cA_{2}=\left\{(x^1,x')\  | \  (x^1,x')\in \big( \frac{t_1+t_2}{2}, +\infty\big)\times \R^{N-1}\setminus B_1\big((t_1+t_2)e_1\big) \right\}.$$

Observe that
$$\int_{\cA_{1}\cup\cA_{2}}(\frac{1}{|t_{1}e_{1}-y|^{N+2s}}-\frac{1}{|t_{2}e_{1}-y|^{N+2s}})dy=0.$$
Since $|t_{1}e_{1}-y|>|t_{2}e_{1}-y|$ for any $y \in \cA_0$,  then it deduces that
$$\varphi(t_1)-\varphi(t_2)=\int_{\cA_0}(\frac{1}{|t_{1}e_{1}-y|^{N+2s}}-\frac{1}{|t_{2}e_{1}-y|^{N+2s}})dy<0. $$
and then
\begin{align*}
\varphi_{B_1}(x) = \int_{B_1^c}\frac{dy}{|x-y|^{N+2s}}
 &=  \int_{B_1^c(x)}\frac{dz}{|z|^{N+2s}}
\\[2mm]&= \int_{B_{\frac{1}{1-|x|}}^c(\frac{x}{1-|x|})}\frac{(1-|x|)^N d\tilde z}{(1-|x|)^{N+2s} |\tilde z|^{N+2s}}
\\[2mm]&= \frac{1}{(1-|x|)^{2s}}\int_{B_{\frac{1}{1-|x|}}^c(\frac{x}{1-|x|})}\frac{d\tilde z}{|\tilde z|^{N+2s}}.
\end{align*}
Combining with $\bigcap {B_{\frac{1}{1-|x|}}^c (\frac{x}{1-|x|})}=(-\infty, -1)\times \R^{N-1}$
and
\begin{align*}
    \int_{(-\infty, -1)\times \R^{N-1}}\frac{d\tilde z}{|\tilde z|^{N+2s}}
    &= \int_{-\infty}^{-1}{d\tilde z_1}\int_{\R^{N-1}}\frac{d\tilde z'}{(|\tilde z_1 |^2+|\tilde z'|^2 )^{\frac{N+2s}{2}}}
     \\[2mm]&= \int_{-\infty}^{-1}{d\tilde z_1}\int_{\R^{N-1}}\frac{{\tilde z_1}^{N-1}dt'}{|\tilde z_1 |^{N+2s}{(1+|t'|^2 )^{\frac{N+2s}{2}}}}
     \\[2mm]&= \int_{\R^{N-1}}\frac{dt'}{{(1+|t'|^2 )^{\frac{N+2s}{2}}}}\int_{-\infty}^{-1}\frac{d\tilde z_1}{{\tilde z_1}^{2s+1}} \, =: c_1,
\end{align*}
it deduces  (\ref{ext 1-beh}).
\hfill$\Box$\medskip

\subsection{Potential inequalities}

For $r_0>0$, denote $\Phi_{r_0}$  the Green kernel of $(-\Delta)^s$ in $B_{r_0}$ with the zero Dirichlet boundary condition in $\R^N\times\R^N\setminus (B_{r_0}\times B_{r_0})$, observe that
\begin{equation}\label{4.2-ap}
\Phi_{r_0}(x,y)\leq c_{31} |x-y|^{2s-N}
\end{equation}
for some $c_{31}>0$ independent of $r_0$.

\begin{lemma}\label{embedding-a}
Assume that $s\in(0,1)$  and integer $N\geq 2$.

$(i)$ If $$\frac1q<\frac{2s}N,$$ then there exists some $c_{32}>0$
such that
\begin{equation}\label{4.1-ap}
\|\Phi_{r_0}\ast h\|_{\infty}\le c_{32}\|h\|_{q};
\end{equation}

$(ii)$ If $$\frac1q\le \frac1r+\frac{2s}N,\quad
\quad q>1,$$
   then
 there exists some $c_{33}>0$ such that
\begin{equation}\label{4.2-a}
\|\Phi_{r_0}\ast h\|_{r}\le c_{33}\|h\|_{q}.
\end{equation}

$(iii)$ If $$1<\frac1r+\frac{2s}N,$$
   then
 there exists some $c_{34}>0$ such that
\begin{equation}\label{4.02-a}
\|\Phi_{r_0}\ast h\|_{r}\le c_{34}\|h\|_1.
\end{equation}
\end{lemma}
\noindent{\bf Proof.} Together with (\ref{4.2-ap}), we  apply  Hardy-Littlewood-Sobolev theorem for the fractional integration \cite[Chapter 5, section 1]{Stein}. For the convenience of the readers, we provides the details of the proof. \smallskip

    {\it Proof of (\ref{4.1-ap})}.
For any $x\in\Omega$ and $q'=\frac q{q-1}$, by
H\"older inequality and (\ref{4.2-ap}), it holds that
\begin{eqnarray*}
 \| \Phi_{r_0}  \ast h\|_{\infty}
 &\le& \Big\|\Big(\int_{B_{r_0}} \Phi_{r_0}^{q'}dy\Big)^{\frac1{q'}} \Big(\int_{B_{r_0}(x)}|h(y)|^q dy\Big)^\frac1q\Big\|_{L^\infty(\Omega)}
 \\[2mm]&\le&c_{35}\|h\|_{q}\Big(\int_{B_{r_0}(x)} \frac1{|x-y|^{(N-2s)q'}}dy\Big)^\frac1q
 \\[2mm]&\leq& c_{36}\|h\|_{q},
\end{eqnarray*}
by the fact that  $$\frac1q<\frac{2s}N,\qquad
  (N-2s)q'<N$$ and
$$\int_{B_{r_0}(x)} \frac1{|x-y|^{(N-2s)q'}}dy<+\infty.$$\smallskip

{\it Proof of (\ref{4.2-a}) and (\ref{4.02-a}) with $r\le q$.}  By  Minkowski inequality, we have that
\begin{eqnarray*}
\|(\Phi_{r_0} \ast h\|_{r}&=  &\|h\ast  \Phi_{r_0}  \|_{r}
\\[2mm]& \le& c_{37}\Big[\int_{\R^N}\Big(\int_{B_{r_0}}\frac{|h|(x-y)\chi_{_{B_2(0)}}(y)}{|y|^{N-2s}}dx\Big)^rdy\Big]^{\frac1r}
  \\[2mm]& \le&
c_{38}\Big[\int_{\R^N}\int_{\R^N}\frac{|h(x-y)|^r  \chi_{_{B_2(0)}}(y)}{|y|^{(N-2s)r}}dxdy\Big]^{\frac1r}
 \\[2mm]& \le& c_{39}\Big[\int_{\R^N}\int_{\R^N} |h(x-y)|^r dx \frac{\chi_{_{B_2(0)}}(y)}{|y|^{N-2s}}dy\Big]^{\frac1r}
\\[2mm]& \le&  c_{40}\|h\|_{L^r(\R^N)}.
\end{eqnarray*}

{\it Proof of (\ref{4.2-a}) and (\ref{4.02-a}) with
$r> q\ge1$ and $\frac1q\le \frac1r+\frac{2s}N$}. We claim that
if $r>s$ and $\frac1{r^*}=\frac1{q}-\frac{2s}N$, the mapping
$h\to (\Phi_{r_0} \eta_0)\ast h$ is weak-type $(q, r^*)$ in the sense that
\begin{equation}\label{weak rs}
\Big|\{x\in\R^N: |(\Phi_{r_0} \eta_0)\ast h |>t\}\Big|\le
\Big(A_{q,r^*}\frac{\|h\|_{L^{q}(\Omega)}}{t}\Big)^{r^*},\qquad h\in
L^q(B_{r_0}),
\end{equation}
for all $t>0$, where  $A_{q,r^*}>0$.

For $\nu>0$, we denote
$$
G_1= \Phi_{r_0}\eta_0 \chi_{_{B_\nu }},\quad \quad G_2=\Phi_{r_0}\eta_0 \chi_{_{B_\nu ^c}}.
$$
it deduces that
\begin{eqnarray*}
&&\Big|\{x\in B_{r_0}: |(\Phi_{r_0} \eta_0)\ast h(x)|>2t\}\Big|
\\[2mm]&\le & \Big|\{x\in\R^N: | G_1\ast h (x)|>t\}\Big|+\Big|\{x\in\R^N: |G_2\ast h(x)|>t\}\Big|.
\end{eqnarray*}

One hand, by  Minkowski inequality, we have that
\begin{eqnarray*}
\Big|\{x\in\R^N: |G_1\ast h(x)|>t\}\Big|&\le&
\frac{\|G_1\ast h\|^s_{s}}{t^s}=\frac{\|h\ast G_1\|^s_{s}}{t^s}
\\[2mm]&\le&\frac{[\int_{\R^N}(\int_{\R^N} |h(x-y)|^sdx)^{\frac1s}|y|^{2\alpha-N}\chi_{_{B_\nu}}(y)dy]^s}{t^s}
\\[2mm]&\le&\frac{\|h\|^s_{s} }{t^s}\int_{B_\nu}|y|^{2\alpha-N}dy
=c_{41}\nu^{2\alpha}\frac{\|h\|^s_{s} }{t^s}.
\end{eqnarray*}
 On the other hand, direct computation shows that
\begin{eqnarray*}
\|G_2\ast h\|_{\infty}&\le&c_{42} \Big\|\int_{\R^N}
\chi_{_{B_\nu^c}}(x-y)(\Phi_{r_0}\eta_0)|h(y)|dy\Big\|_{ \infty}
\\[2mm]&\le& \Big(\int_{\R^N} |h(y)|^qdy\Big)^{\frac1q} \Big\|\Big(\int_{B_2(x)}\chi_{_{B_\nu^c}}(\Phi_{r_0}\eta_0)^{q'}  dy\Big)^{\frac1{q'}}\Big\|_{\infty}
\\[2mm]&\le& \|h\|_{q}\| \Phi_{r_0}\eta_0 \chi_{_{B_\nu^c}}\|_{q'},
\end{eqnarray*}
where $q'=\frac q{q-1}$ if $q>1$, if not, $q'=\infty$.

Since $$\|\Phi_{r_0}\eta_0 \chi_{_{B_\nu^c}}\|_{L^{q'}(\R^N)}
=\Big(\int_{B_2\setminus B_\nu}|x|^{(2s-N)q'}dx\Big)^{\frac1{q'}}
=c_{43}\nu^{2s-\frac Nq},$$
letting $\nu=(\frac t{c_{43}\|h\|_{q}})^{\frac1{2s-\frac Nq}}$, we have that
$$\|G_2\ast h\|_{\infty}\le t,$$
that is,  $$\Big|\{x\in\R^N: |G_2\ast h(x)|>t\}\Big|=0.$$
Then
$$\Big|\{x\in \R^N: | (\Phi_{r_0}\eta_0)\ast h|>2t\}\Big|\le \frac{c_{44}\|h\|^q_{q}\,\nu^{2sq}}{t^q}
\le \frac{ c_{45}{\|h\|_{q}^{r^*}}}{t^{r^*}}.$$
The case  $(ii)$ and $(iii)$ with $r>s$ follows by  Marcinkiewicz Interpolation Theorem.
 \hfill$\Box$\bigskip

    \bigskip

\noindent {\bf  Conflicts of interest:} The authors declare that
they have no conflicts of interest regarding this work.

\medskip

\noindent {\bf Data availability:} This paper has no associated data.

\medskip

\noindent{\bf  Acknowledgements:}  H. Chen is supported by  NNSF of China (Nos. 12071189, 12361043) and Jiangxi Natural Science Foundation (No. 20232ACB201001).

\end{document}